\documentstyle[12pt,amstex,amssymb]{article}

\def\fin{e.\ g.\ } 
\def\this{i.e.\ } 
\def\h{\hbar}
\def\p{\partial}
\def\w{\wedge}

\def\dim{\operatorname{dim}}

\def\Im{\operatorname{Im}}
\def\Re{\operatorname{Re}}
\def\Res{\operatorname{Res}}
\def\mod{\operatorname{mod}}
\def\fgt{\operatorname{fgt}}
\def\ev{\operatorname{ev}}
\def\mk{\operatorname{mk}}

\newcommand{\CC}{{\Bbb C}}
\newcommand{\RR}{{\Bbb R}}

\newcommand{\ZZ}{{\Bbb Z}}
\newcommand{\QQ}{{\Bbb Q}}
\newcommand{\NN}{{\Bbb N}}
\newcommand{\lan}{\langle}
\newcommand{\ran}{\rangle}

\renewcommand{\a}{\alpha}
\renewcommand{\b}{\beta}

\renewcommand{\d}{\delta}

\renewcommand{\l}{\lambda}
\newcommand{\m}{\mu}

\renewcommand{\t}{\tau}

\newcommand{\z}{\zeta}

\newcommand{\D}{\Delta}
\renewcommand{\L}{\Lambda}

\newcommand{\cald}{{\cal D}}

\newcommand{\calg}{{\cal G}}

\newcommand{\cali}{{\cal I}}
\newcommand{\calj}{{\cal J}}

\newcommand{\call}{{\cal L}}

\newcommand{\caln}{{\cal N}}
\newcommand{\calo}{{\cal O}}

\newcommand{\cals}{{\cal S}}
\newcommand{\calt}{{\cal T}}

\newcommand{\calv}{{\cal V}}

\title{A mirror theorem for toric complete intersections}
\author{Alexander Givental \thanks{ Research is partially supported 
by Alfred P. Sloan Foundation and by NSF grant DMS-9321915}}
\date{January 27, 1997 
\abstract{ \footnotesize 
We prove a generalized mirror conjecture for non-negative
complete intersections in symplectic toric manifolds. Namely, we express
solutions of the PDE system describing quantum cohomology of such a manifold
in terms of suitable hypergeometric functions.}}

\begin{document}
\maketitle

{\bf 0. Introduction.}
Let $X$ denote a non-singular compact K\"ahler toric variety with the Picard
number $k$. The variety $X$ 
can be obtained by the symplectic reduction of the standard
Hermitian space $\CC ^N$ by the action of a subtorus $T^k$ in the torus $T^N$
of diagonal unitary matrices on a generic level of the momentum map
(see Section $3$). The coordinate hyperplanes in $\CC ^N$ are $T^N$-invariant,
and their $T^k$-reductions on the same level of the momentum map define
$N$ compact toric hypersurfaces in $X$. We denote $u_1,...,u_N$ the classes
in $H^2(X)$ Poincare-dual to the fundamental cycles of these hypersurfaces.
It is known that $H^2(X)$ is a free abelian group of rank $k$ spanned by
$u_1,...,u_N$, that it multiplicatively generates the ring $H^*(X)$, and 
that the $1$-st Chern class $c(\calt _X)$ of the tangent bundle to $X$ is
equal to $u_1+...+u_N$. 

Let us consider the sum $\calv$ of $l\geq 0$ non-negative line bundles over 
$X$ with the $1$-st Chern
classes $v_1,...,v_l$ and denote $Y$ the non-singular complete intersection
in $X$ of dimension $N-k-l$ defined by global holomorphic sections 
of these line bundles. The inclusion $Y\subset X$ induces the homomorphism
$H^2(X)\to H^2(Y)$       
The cohomology ring $H^*(Y)$ contains the subring 
multiplicatively generated by the classes $u_1,...,u_N$. We denote this
subring $H^*(\calv )$. It carries the
Poincare pairing $\lan f, g\ran =\int _{[Y]} fg =\int _{[X]} fgv_1...v_l$
non-degenerate over $\QQ $. The $1$-st Chern class $c(\calt_Y)$ of the tangent
bundle to $Y$ is equal to $u_1+...+u_N-v_1-...-v_l$.    

\medskip
    
Given a compact K\"ahler manifold $Y$, the Gromov - Witten theory \cite{RT}
associates to it the {\em quantum cohomology algebra} and the {\em quantum
cohomology $\cald $-module}.

Let us introduce the semigroup $\L \subset H_2(Y)$ of fundamental classes of
compact holomorphic curves in $Y$. The quantum cohomology algebra is 
a skew - commutative associative deformation of the cohomology algebra 
$H^*(Y, \QQ [[\L ]])$ with coefficients in a formal completion of the semigroup
algebra $\QQ [\L ]$. The structural constants $\lan a \circ b , c\ran $
of the quantum multiplication $\circ $ are formal series 
$\sum _{d\in \L} n_d(a,b,c) q^d $ where $q^d$ is the element of the semigroup
ring corresponding to the homology class $d\in \L$. The coefficient 
$n_d(a,b,c)$ has the meaning of the number of holomorphic maps $\CC P^1 \to Y$
representing the homology class $d$ and sending $0,1,\infty \in \CC P^1$
to given cycles in $Y$ Poincare-dual to the classes $a,b,c \in H^*(Y)$ 
respectively. We refer to \cite{BM, RT} for variants of actual definitions
which employ {\em intersection theory} in Kontsevich's compactifications of 
moduli spaces of stable maps $\CC P^1\to Y$ (see Section $1$). In particular
$n_d(a,b,c)=0$ unless the total degree of the classes $a,b,c$ equals the
real dimension of the fundamental class in the space of maps 
$\CC P^1 \to Y$ representing the class $d$. This gives rise to the grading
$\deg q^d =2 \int _{[d]} c(\calt_Y)$ in the quantum cohomology algebra.

In addition to the associativity identity the structural constants of the 
quantum multiplication satisfy some integrability condition which can
be formulated as compatibility of certain system of linear differential
equations. Let $(t_1,...,t_k)$ denote coordinated on $H^2(Y)$ with respect
to a basis $(p_1,...,p_k)$ of integral symplectic classes. The basis $q^d$
in the semigroup ring $\ZZ [\L ]$ can be identified with 
$q_1^{d_1}...q_k^{d_k}=\exp (t_1d_1+...+t_kd_k)$ where $(d_1,...,d_k)$ are
coordinates of $d\in H_2(Y)$ with respect to the dual basis. Denote $t_0$ the 
coordinate on $H^0(Y)$. The operators $p_i\circ $ of quantum 
multiplication act on vector functions $s(t_0,t)$
with values in $H^*(Y,\CC )$.
In these notations the integrability condition reads:
\[ \h \frac{\p }{\p t_0} s = s, \  \h \frac{\p}{\p t_i} s = p_i \circ s,
i=1,...,k,\]
form a consistent system of linear PDEs. 

Solutions $s$ to the PDE system can also be interpreted in terms of 
intersection theory on spaces of stable maps (see Section 1). The 
$\cald $-module corresponding to this syslem is generated by a single
formal vector-function $J_Y(t,\h ^{-1})$ with coefficients in $H^*(Y,\QQ )$
(see Section $1$). It has the following property  \cite{Gi3}: whenever
a scalar differntial operator $D( \h \p/\p t, \exp t, \h)$ with coefficients in
$\CC [\h ][[\L ]]$ annihilates 
(all componenets of the vector-function) $J_Y$, the symbol $D(p\circ ,q,0)$
vanishes in the quantum cohomology algebra of $Y$.

A construction of 
$J_Y$ in terms of intersection theory on moduli spaces of stable maps is
given in Section $1$. By the very definition the function $J_Y$ has the 
the asymptotical expansion 
\[ J_Y=e^{(t_0+p_1t_1+...p_kt_k)/\h } (1+o(1/\h))  .\]

\medskip

In applications to a toric complete intersection $Y\subset X$ we will detect 
degrees
of holomorphic curves in $Y$ by their degrees in $X$. This means that $\L $ will 
denote the
semigroup of classes of compact holomorphic curves in $X$, $t_1,...,t_k$ will be 
coordinates
on $H_2(X)$, and $J_Y$ --- a formal series in $q_i = e^{t_i}$. Notice that in 
the case where 
$v_1,...,v_l$ are positive and $\dim Y >2$, the map $H_2(Y)\to H_2(X)$ is an 
isomorphism and thus
the notation $J_Y$ has the same meaning as before. However for semi-positive 
$v_1,...,v_l$ (and 
$\dim Y >1 $) the map $H_2(Y)\to H_2(X)$ can have a non-trivial kernel, and 
therefore $J_Y$ is
obtained from the vector-function discussed in the previous paragraphs by the 
restriction operation.
\footnote{ I would like to thank V. Batyrev and B. Kim who pointed me to this 
subtlety in the 
semi-positive case.}
   
We will denote $J$ the orthogonal projection of the vector-function $J_Y$ to the 
subalgebra 
$H^*(\calv, \QQ ) \subset H^*(Y,\QQ )$. 

\medskip

Introduce the formal vector-function $I(t,\h ^{-1})$ with values in
$H^*(\calv, \QQ ) $:
\[ (*) \ \ \ \ I=e^{(t_0+p_1t_1+...+p_kt_k)/\h }
\sum_{d\in \L} e^{t_1d_1+...+t_kd_k} \ \times \]
\[\ \times \ \frac{\Pi_{a=1}^l \Pi_{m=-\infty }^{L_a(d)} (v_a+m\h )
\ \Pi_{j =1}^N \Pi_{m=-\infty}^0 (u_j+m\h ) }
{\Pi_{a=1}^l \Pi_{m=-\infty}^0 (v_a+m\h) 
\ \Pi_{\a=1}^N \Pi_{m=-\infty }^{D_j(d)} (u_j+m\h )} ,\]
where $D_j(d)= \int_{[d]} u_j$, $L_a(d)= \int _{[d]} v_a$.
The function $I$ is the product of $e^{(t_0+p\log q)/\h }$ with a
power series in $q=e^t$ supported in $\L $ and weighted-homogeneous of
degree $0$ with respect to the grading 
\[ \deg \h =\deg u_j =\deg v_a=1, \deg q^d = \int_{[d]} c(\calt_Y) .\]
The function $J$ has the same properties (see Section $1$).

\medskip

{\bf Theorem 0.1.} {\em Let $Y$ be a non-singular toric complete intersection
with the non-negative $1$-st Chern class $c(\calt_Y)$. Then the formal 
vector-functions $I$ and $J$ with coefficients in $H^*(\calv ,\QQ )$ coincide
up to a triangular weighted-homogeneous change of variables:
\[t_0 \mapsto t_0 + f_0(q) \h + h(q),\ \log q_i\mapsto \log q_i + f_i(q) \]
where $h, f_0, f_1,..., f_k $ are weighted-homogeneous power series 
supported in $\L - 0$ and $\deg f_0 =\deg f_i = 0$.}

\medskip

{\em Remarks.} $1)$ 
Comparisson of the asymptotical expansions for the components
in $H^0(\calv, \QQ ) $ and $H^2(\calv ,\QQ) $ of the functions $I$ and $J$ 
in orders $h^0$ and $h^{-1}$
uniquely determines the change of variables that transforms $I$ to $J$
(see Section $7$). 

 
$2)$ The assumption $c(\calt _Y) \geq 0$ guarantees that $\deg q^d \geq 0$ for
all $d\in \L-0$. If $c(\calt _Y) $ is positive, then the change of
coordinates transforming $I$ to $J$ reduces to the multiplication by
$\exp (h(q)/\h)$ where $h(q)$ is a polynomial supported at 
$\{ d\in \L | \int_{[d]} c(\calt_Y) =1 \} $. Thus the functions $I$ and $J$
coincide whenever this support is empty (for instance, if $c(\calt_Y)$ is 
a multiple of a positive integer class).  

$3)$ In the opposite case $\sum v_a =\sum u_j$ the manifold $Y$ is a
{\em Calabi-Yau} toric complete intersection 
(in the broad sense that $c(\calt_Y)=0$). 
Then $\deg q_i=0$
for all $i$, $h(q)=0$ and
\[ \exp f_0 =\sum _{d: D_j(d)\geq 0} q^d 
\frac{L_1(d)!...L_l(d)!}{D_1(d)!...D_N(d)!} .\]
This and other components of the vector-function $I$ are generalized 
hypergeometric functions.

\medskip

{\bf Example: toric manifolds themselves.} Taking $l=0$ in Theorem $0.1$ we
arrive at the case $Y=X$ of Gromov -- Witten theory on non-singular toric
symplectic manifolds with $c(\calt_X)\geq 0$. We describe further 
simplifications which occur in this case.

(a) If the anti-canonical class of the toric manifold is positive then
$I=J_X$. The proof consists in elementary verification of the asymptotocal
expansion $I=\exp((t_0+p\log q)/\h)\ (1+o(1/\h ))$. In particular, symbols of
differential operators annihilating $I$ yield the relations 
\[ u_1^{D_1(d)}...u_N^{D_N(d)}=q^d \]
in the quantum cohomology algebra of $X$ (see Corollary $0.4$ below).
This result agrees with the multiplicative structure in Floer cohomology
of the loop space $LX$ constructed in \cite{Gi} by the method of 
generating functions (and based on discretization of loops) and with
the results in \cite{B1, MP} on quantum cohomology of positive toric
manifolds.

(b) If the anti-canonical class of the toric manifold is semi-positive then
\[ I=e^{(t_0 +p\log q)/\h} \ [1+p_1f_1(q)/\h+...+p_kf_k(q)/\h \ +\  o(1/\h)]\] 
and differs from $J_X$ only by the change of variables 
$\log q_i\mapsto \log q_i + f_i(q) $. We illustrate this case by the
following instructive example which also exhibits {\em symplectic} invariance
of Gromov -- Witten theory.

Let $X_1$ and $X_2$ be the toric $3$-folds obtained by projectivization
of the following $3$-dimensional vector bundles over $\CC P^1$:
$\calo (-1)\oplus \calo (-1) \oplus \calo $ for $X_1$ and
$\calo (-2) \oplus \calo \oplus \calo $ for $X_2$. The bundles are
topologically (and symplecticly) equivalent and the manifolds $X_1$ and $X_2$
are symplectomorphic. They are not isomorphic however as complex manifolds,
and we shell see how the same Gromov -- Witten invariants emerge from
different computational procedures.

The manifolds $X_1$, $X_2$ are obtained by factorization of $\CC^5$ by 
the torus $T_{\CC }^2$ embedded into the $5$-torus of diagonal matrices
as prescribed by the matrices $M^t_1$ and $M^t_2$ respectively:
\[ M_1 =\left[ \begin{array}{rrrrr} 
1 & 1 & 0 & -1 & -1 \\
0 & 0 & 1 &  1 &  1 \end{array} \right] , \ 
   M_2=\left[ \begin{array}{rrrrr} 
1 & 1 & 0 & 0 & -2 \\
0 & 0 & 1 & 1 &  1 \end{array} \right] .\]
The columns of the matrix $M$ represent the classes $u_1,...,u_5$ as linear
combinations of a basis $p_1, p_2$ in $H^2(X)$: for $X_1$ we have
$u_1=u_2=p_1, u_3=p_2, u_4=u_5=p_2-p_1$. Using upper-case notations
for $X_2$ we get $U_1=U_2=p_1, U_3=U_4=p_2, U_5=p_2-2p_1$. We refer the reader
to Section $3$ for a detailed combinatorial description of cohomological
invariants of toric manifolds. Using that description we find the 
multiplicative relations $u_1u_2=0, u_3u_4u_5=0$ and $U_1U_2=0, U_3U_4U_5=0$ 
in the cohomology algebra and observe that they are effectively the same:
\[ H^*(X)=\ZZ [p_1,p_2]/(p_1^2, p_2^3-2p_2^2p_1) .\]
The K\"ahler cone in both cases consists of $p_1t_1+p_2t_2$ with $t_1,t_2>0$.
The anti-canonical class $\sum u_j =\sum U_j = 3p_2$, 
and the semigroup algebra $\QQ[[\L]]$ of identifies with $\QQ [[q_1, q_2]]$ 
with the grading $\deg q_1=0, \deg q_2 = 3$.  

The series $I$ corresponding to $X_1$ has the form
\[ I_1= e^{(t_0+p_1\log q_1 + p_2\log q_2)/\h} \times \]
\[ \sum_{d_1,d_2=0}^{\infty}
\frac{q_1^{d_1}q_2^{d_2}\ \Pi_{m=-\infty}^0 (p_2-p_1+m\h)^2}
{\Pi_{m=1}^{d_1} (p_1+m\h)^2 \ \Pi_{m=1}^{d_2} (p_2+m\h) \ 
\Pi_{m=-\infty }^{d_2-d_1} (p_2-p_1+m\h)^2} \]
It has the asymptotics $\exp((t_0+p\log q )/\h)\ [ 1+o(1/\h )]$
and thus coincides with $J_X$. 

The series $I_1$ is annihilated by the differential operators
\[ \D _1=(\h \frac{\p}{\p \log q_1})^2 - q_1 (\h \frac{\p}{\p \log q_2} -
\h \frac{\p}{\p \log q_1})^2, \] 
\[ \D _2=\h \frac{\p }{\p \log q_2} (\h \frac{\p}{\p \log q_2} -
\h \frac{\p}{\p \log q_1})^2 - q_2 .\]
Thus the rlations 
\[ p_1^2=q_1 (p_2-p_1)^2, \ p_2 (p_2-p_1)^2 = q_2 \]
describe the quantum deformation of the cohomology algebra of $X$.

The series $I$ corresponding to $X_2$ has the form (we use the
upper-case notations $Q_1,Q_2$ instead of $q_1,q_2$):
\[ I_2 =e^{(t_0+p_1\log Q_1 +p_2\log Q_2)/\h} \times \]
\[ \sum_{d_1,d_2=0}^{\infty }
\frac{Q_1^{d_1}Q_2^{d_2}\ \Pi_{m=-\infty}^0 (p_2-2p_1+m\h)}
{\Pi_{m=1}^{d_1} (p_1+m\h)^2\ \Pi_{m=1}^{d_2} (p_2+m\h)^2\ 
\Pi_{m=-\infty}^{d_2-2d_1} (p_2-2p_1+m\h)} .\]
It has the asymptotics $\exp((t_0+p\log Q)/\h)\ 
[1 + (2p_1-p_2) f(Q_1)/\h  + o(1/\h )] $ where
\[ f(Q_1)=\sum_{d_1=1}^{\infty } \frac{(2d_1-1)!}{(d_1!)^2} Q_1^{d_1} .\]
We have to put therefore $q_1=Q_1\exp(2f(Q_1)), \ q_2=Q_2\exp(-f(Q_1))$.

In fact the inverse change of variables is given by the simple formulas
\[ Q_1=\frac{q_1}{(1+q_1)^2}, \ Q_2=q_2(1+q_1), \]
 \[ Q_2 \frac{\p}{\p Q_2}= q_2 \frac{\p}{\p q_2},\ 
Q_2 \frac{\p}{\p Q_2} - 2Q_1\frac{\p}{\p Q_1} =
\frac{1+q_1}{1-q_1} (q_2\frac{\p }{\p q_2}-2q_1\frac{\p }{\p q_1}) .\]
With these formulas at hands it is straightforward to check that 
the change of variables transforms the system of differential equations
\[ \h ^3 (Q_2\frac{\p}{\p Q_2})^2 
(Q_2\frac{\p}{\p Q_2} - 2Q_1\frac{\p}{\p Q_1})\ I=Q_2\ I \]
\[ (Q_1 \frac{\p}{\p Q_1})^2 \ I = 
Q_1 (Q_2\frac{\p}{\p Q_2}-2Q_1\frac{\p}{\p Q_1})
(Q_2\frac{\p}{\p Q_2} -2Q_1\frac{\p}{\p Q_1} -1)\ I \]
satisfied by $I=I_2$ to the system $\D_1 I = 0, \D_2 I=0$
satisfied by $I=I_1$. This guarantees that $I_2(Q(q))=I_1(q)$ and finally
gives rise to the same description of the quantum cohomology algebra.
$\square $                 

\medskip

Our proof of Theorem $0.1$ is based on an equivariant generalization of the
Gromov -- Witten theory. 

Let $\calv $ be a holomorphic $l$ - dimensional
vector bundle over the $n$ - dimensional
K\"ahler manifold $X$ equivariant with repect to a hamiltonian Killing
action of a torus $G$. It was explained in \cite{Gi3} how to extend the 
Gromov -- Witten theory to the {\em $G$ - super - manifolds} $(X,\calv)$
of dimension $(n,l)$ 
in the case of {\em convex} $\calv $ (\this bundles with all fibers  
spanned by global holomorphic sections). 
\footnote{The idea of the
construction is due to M. Kontsevich (see \cite{Kn,Gi3})
while the terminology of super-manifolds in this context was introduced
by A. Schwarz \cite{Sch}.} 

The equivariant quantum
cohomology algebra of the super - manifold $(X,\calv )$ is a deformation of the
cup - product in the equivariant cohomology algebra $H^*_G(X,\QQ)$ provided
with the Poincare pairing $\lan f ,g\ran = \int_{[X]} fg\ Euler(\calv )$
with values in $H^*_G(point)=H^*(BG,\QQ)$. Here the equivariant top Chern class
$Euler (\calv)$ is assumed to be invertible over the field of fractions of the
polynomial algebra $H^*(BG,\QQ)\simeq \QQ [\l_1,...,\l_{\dim G}]$. 
Respectively, the equivariant quantum 
cohomology algebra and the quantum cohomology $\cald $-module 
of the super-manifold are defined over this field of fractions. It is important
however that the structural constants $\lan a \circ b, c\ran $ and the 
components $\lan J_{\calv} , a \ran $ of the corresponding solution
vector-function are defined over the polynomial algebra $\CC [\l ]$
(as some equivariant Poincare pairings with the top Chern classes of 
suitable vector bundles over the moduli spaces of stable maps to $X$,
see Section $1$). In the non-equivariant limit to $\l =0$ the algebra
$(H^*_G(X)/ ker \lan \cdot ,\cdot \ran )$ degenerates to $H^*(\calv )$,
and the corresponding structural constants and solutions turn into  
their counterparts in the Gromov-Witten theory on the complete intersection 
$Y\subset X$ defined by a generic holomorphic section of $\calv $.
This allows to obtain Theorem $0.1$ as the specialization to $\l =0$ of 
the following result about toric super-manifolds. 

\medskip

Let us consider the equivariant cohomology
algebra $H_{T^N}^*(X)$ of the toric manifold $X=\CC ^N// T^k$
provided with the action of the torus $T^N$ of diagonal unitary
matrices. The coefficient algebra $H^*(BT^N)$ of the equivariant theory
is canonically isomorphic to the algebra $\ZZ [\l_1,...,\l_N ]$ of 
polynomial functions on $Lie T^N$. The algebra $H^*_{T^N}(X)$ 
is multiplicatively
generated over $\ZZ [\l ]$ by the equivariant counterpart
$p_1,...,p_k \in H^2(X)$ of a basis in $H^2(X)$. The equivariant classes
$u_1,...,u_N \in H^2_G(X)$ Poincare-dual to the $G$-invariant toric
hypersurfaces corresponding to the coordinate hyperplanes in $\CC ^N$
are some linear combinations 
\[ u_j=\sum _{i=1}^k p_i m_{ij} - \l _j, \ j=1,...,N, \] 
of these generators. The relations between the generators $p_i$ can be
written in the form (see Section $3$): 
\[ u_{j_1}...u_{j_r}=0, \ j_1<...<j_r ,\] 
whenever the toric hypersurfaces corresponding to $u_{j_1},..., u_{j_r}$ have
empty intersection.

Let $\calv $ denote the direct sum of $l\geq 0$ non-negative $T^N$-equivariant
line bundles over $X$ with the {\em non-zero} equivariant $1$-st Chern classes
\[ v_a=\sum _{j=1}^N p_i l_{ia} -\l'_a ,\ a=1,...,l .\] 
The variables $(\l',...,\l'_l)$ here stand for the generators of
$H^*(BT^l)$ where the torus $T^l$ acts fiberwise on the bundle $\calv $.
(One may think of the toric super-manifold $(X,\calv )$ as of the symplectic
reduction of the super-space $(\CC ^N, \CC ^l)$ by the torus $T^k$ embedded
into the maximal torus $G=T^N\times T^l$ by means of the matrix 
$(m_{ij} | l_{ia})$.)

Denote $J_{\calv}(t_0, \log q, \h^{-1})$ 
the solution formal vector-function of the
equivariant Gromov -- Witten theory on the super-manifold $(X,\calv)$. 
Define another formal vector-function 
$I_{\calv} (t_0, \log q, \h^{-1})$ by the formal series $(*)$.  
Coefficients of both formal functions are $G$-equivariant cohomology
classes of $X$ over $\QQ [\l ,\l']$. In the limit $(\l,\l' )=0$ the functions
$J_{\calv} $ and $I_{\calv}$ yield $J$ and $I$ respectively.

\medskip

{\bf Theorem 0.2.} {\em Suppose that $1$-st Chern class $\sum u_j - \sum v_a$
of the super-manifold $(X,\calv )$ is non-negative. Then $J_{\calv}$ coincides
with $I_{\calv}$ up to a weighted-homogeneous triangular change of variables:
\[ t_0\mapsto t_0+f_0(q)\h + \sum \l_j g_j(q) + h(q),\ 
\log q_i\mapsto \log q_i + f_i(q), \ i=1,...,k \]
where $f_0, f_i, g_j, h$ are weighted-homogeneous formal $q$-series 
with $\deg f_0=\deg f_i=\deg g_j=0, \deg h=1$ supported at $\L -0$. }
 
\medskip

{\bf Corollary 0.3.} {\em Suppose that a linear differential operator \newline
$D(\h \p /\p \log q, q, \l, \h )$ with coefficients in $\CC [\l, \h][[\L ]]$
annihilates the vector-function $I_{\calv}$ transformed to the new variables. 
Then the relation $D(p\circ ,q,\l,0)$ $= 0$ 
holds in the quantum cohomology algebra of the super-manifold
$(X,\calv )$. In particular, in the quantum cohomology algebra
of the complete intersection $Y\subset X$ (with $\dim Y>1$) we have 
$\lan a , D(p\circ , q, 0, 0)\circ b \ran = 0$ for any 
$a, b \in H^*(\calv) \subset H^*(Y)$.}

\medskip

{\em Remarks.} $4$) The change of variables transforming $I_{\calv}$ to 
$J_{\calv}$ is uniquely determined by the asymptotics of $I_{\calv}$ modulo
$\h^{-2}$. 

$5$) The $q^d$-term in 
$\lan a, b \circ p_{i_1}\circ ...\circ p_{i_r}\ran $ has the following
enumerative meaning: it is the (virtual) number of degree $d$ holomorphic maps
$\CC P^1 \to Y$ which send a given generic configuration of $r+2$ 
distinct points in $\CC P^1$ to the given generic cycles $a, b$ and $r$ given
generic hypersurfaces Poincare-dual to $p_{i_1},...,p_{i_r}$ respectively.   

$6$) Due to a somewhat subtle relationship 
between the quantum cohomology algebras 
of $Y$ and $(X, \calv)$ it seems to be dangerous to deduce enumerative 
corollaries about $Y$ directly from differential equations for $J$ (instead of
$J_{\calv}$ or $J_Y$). We do not know counter-examples however. In many
cases such corollaries can be justified by means of additional dimensional
or Hodge-theoretic arguments. In particular, if $\oplus H^{r,r}(Y) \subset
H^*(\calv, \CC )$ (for example if $Y\subset X$ is a hypersurface of odd 
dimension), we have $\lan a, J_Y \ran =\lan a, J\ran $ for all 
$a\in H^*(Y)$. In this case $D(p\circ ,q,0)=0$ if the differential operator
$D(\h \p /\p \log q, q, \h)$ annihilates $I$ transformed to the
new variables as described in Theorem $0.1$.  

\medskip

Consider now the differential equations satisfied by $I_{\calv}$.
Put 
\[ \p_j=\sum_i m_{ij} \h \p/\p \log q_i -\l_j ,\ j=1,...,N,\] 
\[ \p '_a=\sum _i l_{ia} \p/\p \log q_i -\l'_a, \ a=1,...,l.\]
For each $d\in \L $ with $D_1(d),...,D_N(d)\geq 0$ we introduce the
differential operator (see \cite{Gi2}):
\[ \D _d=\Pi_j \Pi_{m=0}^{D_j(d)-1} (\p_j-m\h) - 
q^d \Pi_a \Pi_{m=1}^{L_a(d)} (\p'_a +m\h) . \]
It follows easily from the detailed description of the equivariant cohomology
algebra $H^*_G(X)$ given in Section $3$ that $\D_d I_{\calv} =0$.
(In fact $I_{\calv}$ may satisfy some stronger differential equations.) 

Denote $\hat{\D}_d$ the polynomial 
\[ u_1^{D_1(d)}...u_N^{D_N(d)}-q^dv_1^{L_1(d)}...v_l^{L_l(d)} \]
in the quantum cohomology algebra of the non-negative toric complete 
intersection $Y\subset X$. 

\medskip

{\bf Corollary 0.4.} {\em Suppose that for the toric complete intersection 
$Y\subset X$ we have $J_{\calv}=I_{\calv}$ (for instance due to the causes 
described in Remark $2$). Then $\lan a, b\circ \hat{\D}_d \ran =0$ for
any $a, b\in H^*(\calv)\subset H^*(Y)$.}   

\medskip

We would like to emphasise the hypergeometric character of the function
$I_{\calv}$. It is easy to see that the differential operators $\D_d$
annihilate also the following hypergeometric integrals:
\[ \int_{\Gamma \subset E_q} \ e^{(\sum u_j - \sum v_a)/\h }
\ u_1^{\l_1 /\h} ... u_N^{\l_N/\h } \ v_1^{\l'_1/\h }...v_l^{\l'_l/\h } \ 
\times \ \]
\[ \ \times \ \frac{d\log u_1 \w ... \w d \log u_N \w d v_1 \w ... \w d v_l}
{d\log q_1 \w ... \w d\log q_k } \ .\]
Here $\Gamma _q $ is a suitable real $N-k+l$-dimensional (non-compact)
cycle in the complex $N-k+l$-dimensional variety 
\[ E_q = \{ (u,v) | \Pi_{j=1}^N u_j^{m_{ij}}=q_i \Pi_{a=1}^l v_a^{l_{ia}},
i=1,...,k\} \]
provided with the local coefficient system $u^{\l /\h } v^{\l' /\h}$.

\medskip

{\bf Example: a mirror theorem.} 
Partition the variables $u_1,...,u_N$ into $l+1$ groups
and denote $F_a$ the sum of $u_j$ in each group, $a=0,...,l$.
Define the matrix $(l_{ia})$ in such a way that $v_a=F_a(u), a=1,...,l$,
correspond to convex line bundles on $X$. In the limit $\l'=0, \l=0$ 
the above
integral evaluated explicitly over $\RR _+^l$ in $dv$ reduces to
\[ \int _{\Gamma '_q \subset X'_q } 
e^{F_0(u)/\h } \frac{d\log u_1\w ... \w d\log u_N}
{(1-F_1(u))...(1-F_l(u))\ d\log q_1 \w ... \w d\log q_k} \]
where $X'_q=\{ u | \Pi _j u_j^{m_{ij}} = q_i , \ i=1,...,k \} $.
Further reduction by Cauchy residue formula yields
\[ (**) \ \ \ \int _{\gamma_q\subset Y'_q} e^{F_0(u)/\h }
\frac{d\log u_1 \w ... \w d\log u_N}{d F_1 \w ... \w dF_l 
\w d\log q_1 \w ... \w d\log q_k} .\]
Here $\gamma_q $ can be understood as Morse-theoretic cycles of 
the function $\Re F_0 $ restricted to the $N-k-l$-dimensional
manifold
\[ 
 Y'_q=\{ u | \Pi_j u_j^{m_{ij}}=q_i, i=1,...,l,\ F_a(u)=1, a=1,...,l \} \ .\]
   
It is not hard to see that all components of the vector-function $I$ are
described by such integrals with suitable $\gamma_q $.
 \footnote{In general the converse is not true --- different toric manifolds 
 $X=\CC ^N//T^k$ obtained by the reduction on different generic levels of
 the momentum map have different cohomology algebras $H^*(\calv)$ and may 
 give rise to different $I$ whose integral representations differ by the 
 choice of cycles only.} 
This constitutes the content of the mirror symmetry between the toric complete
intersections $Y$ and $Y'$. 

\medskip

Indeed:

\medskip
 
(a)  In the case of a Calabi--Yau toric complete intersection $Y$ we have
 $\sum u_j = \sum v_a $ and thus $F_0(u)=0$. The affine varieties $X'_q$
 in this case can be compactified (see \cite{B}) to the toric varietiy 
 $\hat{X}'$ with the momentum polyhedron {\em polar} to that for $X$. 
 The varieties $Y'_q$ 
 compactify to (singular) Calabi--Yau complete intersections 
 $\hat{Y}'_q \subset \hat{X}'$. The forms $d\log u/dv \w d\log q$ extend
 to holomorphic volume forms $\omega _q$ on desingularizations of 
 $\hat{Y}'_q$. Thus the 
 components of the hypergeometric series $I$ identify with the periods
 $\int _{\gamma _{\log q}} \omega _q$ responsible for variations of 
 complex structures in a family Calabi--Yau manifolds 
 birationally isomorphic to $\hat{Y}'_q$.

 According to V.Batyrev \cite{B} the $q$-family of Calabi--Yau manifolds
 is {\em mirror - symmetric} to the original family of toric complete 
 intersections $Y$, and thus our Theorem $0.1$ confirms the {\em mirror
 conjecture} for Calabi -- 
Yau toric complete intersections formulated in detail
 in \cite{BVS}.

\medskip

(b)  In the lecture \cite{Gi2} we suggested a generalization of the mirror 
 conjecture beyond the class of Calabi -- Yau manifolds.
 \footnote{A similar generalization was recently proposed in \cite{EHX}.} 
  In this generalization  
 the quantum cohomology ${\cal D}$-module of a compact K\"ahler manifold $Y$ 
  should be equivalent, up to some change of variables, 
  to the $\cal D$-module generated by the 
  oscillating integrals $\int e^{f_q/\h } \omega_q $ defined by a suitable
  family $(Y'_q, \omega _q, f_q)$ of (possibly non-compact) 
  complex algebraic manifolds $Y'_q$ (of the same dimension as $Y$) 
  provided with holomorphic volume forms $\omega_q$ and holomorphic phase 
  functions $f_q$. This generalisation was confirmed for toric Fano manifolds
  in \cite{Gi2}, for Fano and Calabi--Yau projective complete intersections   
  in \cite{Gi3} and for flag manifolds of the series $A$ in \cite{Gi4}.
  Thus Theorem $0.1$ along with the integral representation $(**)$ proves
  our generalized mirror conjecture for non-negative toric complete 
  intersections described at the beginning of this example. $\square $

\medskip

The remaining part of the paper contains a proof of Theorem $0.2$. It 
exploits some general properties of equivariant Gromov -- Witten invariants
of super-manifolds $(X,\calv)$ described in \cite{Gi3} in the setting of 
{\em convex} K\"ahler manifolds $X$. 
Foundations of the Gromov -- Witten theory for general $X$ were 
developed recently by several groups of authors (see for instance \cite{RT}).
It appears that these results admit staighforward equivariant generalizations.
However there is no ready reference for such a generalization in the 
literature. In Section $1$ we describe 
the properties of the equivariant Gromov--Witten theory 
we use in this paper in the axiomatic form. 
These properties have been verified in \cite{Gi3} for convex $X$.
Thus the proof of Theorems $0.1$ and $0.2$ given in this paper is complete
in the case of non-negative complete intersections in products of projective
spaces (such products are convex), and in the case of more general toric
manifolds still should be complemented by a verification of the axioms.

\medskip 

{\bf 1. Moduli spaces of stable maps.} 
Let $(C,x)$ be a compact connected complex algebraic curve with at most
double singular points and with $r$ pairwise distinct non-singular
marked points $x=(x_1,...,x_r)$. We will assume that $C$ has the arithmetic
genus $\dim H^1(C, \calo )=0$. Two holomorphic maps $(C,x)\to M, (C',x')\to M$
of two such curves to a complex manifold $Y$ are called {\em equivalent} if
they are identified by a holomorphic isomorphism $(C,x)\to (C',x')$. 
A holomorphic map $(C,x)\to M$ is called {\em stable} (see \cite{Kn})
if it does not admit non-trivial infinitesimal automorphisms.
A stable map may have a non-trivial finite group of discrete automorphisms.

The {\em degree} of a holomorphic map $(C,x)\to M$ is defined as the
homology class $d\in H_2(M)$ it represents. 

Denote $M_{r,d}$ the set of equivalence classes of degree $d$ genus $0$ 
stable maps $(C,x)\to M$ with $r$ marked points. 
For a compact projective variety $M$ the set $M_{r,d}$ has a natural structure
of a compact complex algebraic variety (see \cite{Kn, BM}). We will call the
spaces $M_{r,d}$ the {\em moduli spaces} of stable maps.

If $M$ is a homogeneous K\"ahler space of a compact semi-simple Lie group,
the moduli spaces $M_{n,d}$ are known \cite{BM} to be {\em non-singular 
orbifolds} \this local quotients of non-singular spaces by finite groups.
In particular the spaces have natural fundamental cycles providing Poincare
duality over $\QQ $.
For general $M$ it is still convenient to think of the moduli spaces $M_{n,d}$
as of {\em singular orbifolds} \this local quotients of singular spaces by
finite groups. The finite groups in question are the discrete automorphism 
groups of stable maps, and many local constructions in the moduli spaces become
transparent only after passing to their local ``unquotient'' coverings. 
Vector bundles and their characteristic classes provide important examples 
of the orbifold ideology. 
By a vector bundle over $M_{n,d}$ we mean a sheaf which
is locally identified with the sheaf of
invariant sections of a natural local vector bundle on such coverings.
The structural group $G$ of the bundle acts with at most finite stabilizers
on the total space of corresponding principal bundle $P\to M_{n,d}$. In this
case the rational coefficient equivariant cohomology algebra $H^*_G(P,\QQ )$
naturally identifies with $H^*(M_{n,d},\QQ )$. Then the augmentation
homomorphism $H^*(BG,\QQ )\to H^*_G(P,\QQ )$ induced by the $G$-equivariant
map $P\to pt $ defines characteristic classes of the bundle. Talking about 
vector bundles over the moduli spaces, their Euler and Chern classes we
will always have this orbifold subtlety in mind.  

\medskip

According to \cite{FO, RT} the moduli spaces $M_{n,d}$ can be provided
with virtual fundamental cycles ---  rational
homology classes of $M_{r,d}$ satisfying the axioms \cite{KnM} 
of Gromov - Witten theory. We call the virtual fundamental cycles
{\em Gromov - Witten classes}. 
Both the axioms and their realization by means
of Gromov - Witten classes allow dependence of $M$ on parameters.

Suppose that a compact Lie group $G$ acts on a compact K\"ahler manifold $M$
by hamiltonian
automorphisms of the K\"ahler structure. Then $G$ acts naturally on the
Moduli spaces $M_{r,d}$. Let $B\subset BG$ be a 
finite-dimensional approximation to the classifying space $BG$ of principal
$G$-bundles, and $M_B\to B$ --- the associated $M$-bundle. The bundle can be 
considered as a family of compact K\"ahler manifolds. This allows to construct
parametric Gromov - Witten classes in the moduli spaces $(M_{r,d})_B\to B$.
Exhausting $BG$ by the finite-dimensional approximations one can construct
{\em $G$-equivariant} Gromov-Witten classes $[M_{r,d}]$. The recent progress
\cite{FO, RT} in foundations of Gromov -- Witten theory leaves little doubt
that the ``virtual fundamental class'' approach is consistent with the axioms
\cite{KnM, Gi3} of {\em equivariant} Gromov -- Witten theory. 

\medskip

We describe below a variant of the axioms 
in the form convenient for applications in the present paper.

(1) The equivariant Gromov - Witten class defines an
$H^*(BG, \QQ)$ - linear function $\int _{[M_{r,d}]}: H^*_G(M_{r,d},\QQ )\to
H^*(BG,\QQ)$ of homogeneity degree 
$-2 [(c_1(\calt_M),d) + \dim_{\CC} M + r - 3]$.
The moduli spaces $M_{r,0}$ are isomorphic to $M\times \bar{\frak M} _r$ where
$\bar{\frak M} _r$ is the Deligne - Mumford compactification of the moduli 
space ${\frak M}_r$ of ordered $r$-tuples of distinct points on 
$\CC P^1$. The Gromov - Witten classes $[M _{r,0}]$ are (equivariant) 
fundamental cycles of the manifolds $M_{r,0}$ for $r\geq 3$ and are not defined
for $r=0,1,2$. 

{\em Remark.} It is crucial for applications in this paper that the
equivariant Gromov -- Witten class is indeed defined over the polynomial
algebra $H^*(BG,\QQ )$ and not over its field of fractions.

(2) There exist natural equivariant 
{\em forgetting maps} $\fgt_s: M_{r+1,d}\to M_{r,d}$
defined (see \cite{Kn, BM}) by forgetting the marked point 
$x_s, 1\leq s\leq r+1$. 
The fiber of $\fgt_s$ over the point represented by a stable
map $(C,x)\to M$ is canonically isomorphic to the quotient of $C$
by the finite group of automorphisms of the stable map. 

(3) There exist natural equivariant {\em evaluation maps}
$\ev_s:M_{r,d}\to M$ defined by evaluation of a stable map at the marked
point $x_s$. The {\em $3$-point correlators}  
\[ \lan a,b,c \ran := \sum_d q^d \int_{[M_{3,d}]} 
\ev_1^*(a) \ev_2^*(b) \ev_3^*(c)\]
are structural constants $\lan a\circ b, c\ran $ of a (super) commutative 
associative Frobenius algebra structure on $(H^*_G(M), \lan \cdot ,\cdot \ran)$
with the unity $1$. Here $q^d$ stands for the element in the group ring of the
lattice $H_2(M)$ corresponding to $d\in H_2(M)$, and the (super) symmetric 
bilinear form $\lan a , b \ran $ is the (equivariant) Poincare pairing 
$\int_M ab $. The axiom means that the new
$H^*(BG,\QQ [[q]])$-bilinear multiplication $\circ $ on 
$H^*_G(M,\QQ [[q]])$ defined by 
\[ \lan a\circ b, c\ran := \lan a,b,c \ran \ \forall a,b,c \in H^*_G(M) \]
is associative, and that the class $1\in H^0_G(M)$ plays the role of the 
unity: 
\[ \lan a,b,1\ran =\lan a,b\ran \ \forall a,b\in H^*_G(M) .\]
Here $\QQ [[q]]$ is the formally completed semigroup algebra of the
semigroup of degrees $d\in H_2(M)$ of compact holomorphic curves in $M$.
The algebra $(H^*_G(M,\QQ[[q]]), \ \circ )$ is called {\em the (equivariant)
quantum cohomology algebra of $M$}. 
Reduced modulo the maximal ideal in $\QQ[[q]]$, it identifies with 
the ``classical'' cohomology algebra $H^*_G(M,\QQ )$
since $M_{3,0}=M$.

(4) Let $d=(d_1,...,d_k)$ denote coordinates of $d$ with respect to a basis
in $H_2(M)$, and $p_1,...,p_k \in H^2_G(M)$ represent the dual basis in
$H^2(M)/H^2(BG)$. The $H^*(BG)$ - linear operators $p_i\circ $ of quantum 
multiplication by $p_i$ satisfy 
\[ q_i\p (p_j\circ )/\p q_i =
q_j \p (p_i\circ )\p q_j .\] 
These identities along with associativity and
commutativity of the quantum multiplication mean that the following linear
system of parial differential equations is consitent for any value
of the parameter $\h \neq 0 $:
\[ \h \frac{\p }{\p t_0} s = s,\ \ \ \h q_i \frac{\p}{\p q_i} s = p_i\circ s,
\ \ i=1,...,k . \]

(5) When $a$ runs an $H^*(BG)$-basis in $H^*_G(M)$, the following formal 
vector - functions $s_a$ run a basis of solutions to the above PDE system
(see \fin \cite{Gi3}). Define the (equivariant) class $s_a$ of $M$ 
(with appropriate coefficients) by
\[ \forall b\in H^*_G(M) \ \lan s_a, b\ran =
\int_M e^{(t_0+p_1\log q_1+...+p_k\log q_k)/\h } ab + \]
\[ \sum _{d\neq 0} q^d\int_{[M_{2,d}]} 
\frac{\ev_1^*(e^{(t_0+p_1\log q_1+...+p_k\log q_k)/\h } a)}{\h - c} \ 
\ev_2^*(b) .\]
Here $c$ is the (equivariant) $1$-st Chern class of the following line bundle
over $M_{2,d}$ called {\em the universal cotangent line at the $1$-st 
marked point}. The forgetful map $\fgt_3:M_{3,d}\to M_{2,d}$ 
(defined for $d\neq 0$) has the section $\mk_1: M_{2,d}\to M_{3,d}$ defined by
marked point $x_1$ in each fiber. The universal cotangent line is defined as
the conormal bundle to the hypersurface $\mk_1(M_{2,d})$ in $M_{3,d}$. 
The fiber
of this bundle at the point represented by the stable map $(C,x_1,x_2)\to M$
is the cotangent line $T^*_{x_1} C$.

In particular, this axiom implies (see \cite{Gi3}) 
that whenever a differential operator
\[ D(\h q_1\p /\p q_1,...,\h q_k\p /\p q_k, q_1,...,q_k, \h ) \] 
with coefficients in $H^*(BG, \CC [[q ,\h ]]) $ annihilates simultaneously
all functions $\lan s_a, 1\ran , \ a\in H^*_G(M)$, the relation
$D(p_1\circ ,...,p_k\circ , q_1,...,q_k, 0)=0$ defined by the symbol of this
operator holds true in the quantum cohomology algebra of $M$.
The vector-function $J_M$ mentioned in the introduction is defined by
$\forall a \ \ \lan s_a , 1\ran = \lan a, J_M \ran $.    

(6) A holomorphic vector bundle $\calv $ over $M$ is called {\em convex}
if each fiber of $\calv $ is generated by global holomorphic sections of 
$\calv $. The Gromov -- Witten theory of the {\em super-manifold} $(M,\calv )$
is constructed by taking $f\mapsto \int_{[M_{r,d}]} f\ Euler(\calv _{n,d})$
on the role of Gromov -- Witten classes. Here $\calv_{n,d}$ is the vector
bundle
\[ \calv_{r,d} =( \fgt_{r+1} )_* (\ev_{r+1}^* (\calv )) .\]
over the moduli space $M_{n,d}$. 
The fiber of $\calv_{r,d}$ over the point in $M_{r,d}$ represented by
the degree $d$ stable map $\phi: (C,x) \to M$ is $H^0(C,\phi^*(\calv ))$ 
and has the complex dimension $(c_1(\calv ), d)$ since 
$H^1(C,\phi^*(\calv))=0$ due to the convexity of the bundle. 
Evaluation at a marked point defines an epimorphism $\calv_{n,d}\to \calv $
of the bundles over $\ev_s : M_{n,d}\to M$. This allows to introduce 
the correlators $\lan a, b, c\ran , \ \lan s_a, b\ran $ for $a,b,c $
from the quotient $H^*(\calv ) $ of $H^*(X,\QQ )$ by the kernel of the
intersection form $\lan a, b\ran = \int _M a\ b\ Euler(\calv )$ and thus
to define for the super-manifold $(X,\calv )$ the quantum cohomology algebra 
and the quantum cohomology $\cald $-module satisfying the analogues of the
axioms (3),(4),(5) of the super-manifold $(X,\calv )$. Analogous axioms 
hold true in the equivariant setting provided that the convex bundle 
$\calv $ is $G$-equivariant.
Then $Euler (\cdot )$ means the equivariant Euler class. In particular,
the vector-function $J_\calv $ from Introduction is the counterpart of
$J_M$ in the case of the (equivariant) super-manifolds $(X,\calv )$.

(7) Let $Y$ be a non-singular submanifold in $M$ defined as the zero locus of
a section of the convex vector bundle $\calv $. The inclusion $i: Y\subset M$
identifies the Frobenius algebra $(H^*(\calv ), a\mapsto \int_M a\ 
Euler(\calv ))$ with the image of the homomorphism
$i^*: H^*(M,\QQ )\to H^*(Y,\QQ )$. Respectively, the inclusion of the moduli
spaces $Y_{n,d} \to M_{n,i_*d}$ maps the Gromov -- Witten class $[Y_{n,d}]$
to the Gromov -- Witten class corresponding to the super-manifold $(X,\calv )$.
This implies that the correlators $\lan a, b, c\ran , \lan s_a, b\ran $ for $Y$
between the classes $a,b,c$ induced from $M$ can be computed by integration
over the Gromov -- Witten classes $[M_{n,d}]$ against the Euler classes 
of $\calv _{n,d}$. In particular the orthogonal projection $J$ of the
vector-function $J_Y$ to $i^*(H^*(M,\QQ ))$  identifies with (the 
non-equivariant version of) $J_{\calv}$.

(8) The Borel fixed point localization
technique applies to the Gromov - Witten classes $[M_{r,d}]$. 
Consider the action of the torus $G$ on $M_{r,d}$. A fixed point of
this action is represented by a map $\phi: (C,x)\to M$ such that 
each irreducible component of $C$ is mapped either to the fixed point
manifold $M^G$ or onto a complex $1$-dimensional orbit of $G_{\CC }$. 
We will formulate the axiom in the special case where both the set of 
fixed points $M^G$ and the set of $1$-dimensional orbits are finite.
If it is the case,  
each connected component of the fixed point set in
$M_{r,d}$ is (the quotient by a finite group of) a product 
$\Pi \ \bar{\frak M}_n $ of Deligne-Mumford 
configuration spaces $\bar{\frak M}_{n}$. Each factor in this product 
parametrizes configurations of marked and singular points on a
$1$-dimensional connected component of $\phi ^{-1}(M^G)\subset C$. 
On each of the remaining irreducible components of $C$ the map $\phi $  
covers a compactified $1$-dimensional orbit of $G$ (with certain multiplicity
in which case it ramifies at the two compactifying fixed points). We will call
these components {\em edges} of the curve $C$. 

The localization of the equivariant Gromov - Witten class $[M_{r,d}]$ 
at the fixed point component in $M_{r,d}$ is
\[ f\mapsto \frac{1}{|Aut |} \int _{\Pi \ \bar{\frak M}_{n} } 
\frac{f}{\text{Virtual Normal Euler Class}}\ .\]
Here $Aut $ is the automorphism group of a typical stable map $\phi $ from
the fixed point component in question.
 
In order to describe the Virtual Normal Euler Class consider the vector spaces
\[ \caln_i=H^i(C,\phi ^*(\calt_M)), i=0,1,\ \caln_2=H^0(C, \calt_C [x]),\ 
\caln_3= \oplus_{y} T'_y\otimes T''_y\ ,\]
where $y$ runs double points of $C$ situated on the edges, and $T'_y ,T''_y$ 
are the tangent spaces to the two components of $C$ intersecting at $y$.
Dimensions of the spaces $\caln_i$ do not change along the connected fixed
point component of $M_{r,d}$ and thus they form vector bundles over 
$\Pi \ \bar{\frak M}_n$. These bundles carry natural infinitesimal actions of 
$G$. Indeed, this is obvoius for $\caln_0 $ and $\caln_1$. 
The differential $\phi_*$
identifies the space $\caln_2$ of infinitesimal automorphisms of $(C,x)$ 
with a $G$-invariant subspace in $\caln_0$ due to stability of $\phi $. 
For a double point $y$ on the edge $C'$ where $\phi $ is an $m$-fold
covering map, the space $(T_yC)^{\otimes m}$ is identified with the tangent
line to the closed $1$-dimensional $G$-orbit at a fixed point and thus inherits
an infinitesimal action of $G$. The Virtual Normal Euler Class is defined as
the following combination of equivariant Euler classes of the bundles 
$\caln_i$:
\[ \frac{Euler(\caln_0) Euler(\caln_3)}{Euler(\caln_1) Euler(\caln_2)} .\]

{\em Remark.} One can use this formula in order to {\em define} the Gromov
- Witten cycle in the equivariant theory. However the verification of the other
axioms and especially --- the {\em polynomiality} property (1) of the 
Gromov - Witten class to be defined over $H^*(BG,\QQ )=\QQ [\l ]$ ---
becomes then a nontrivial combinatorial task.  

\medskip

{\bf 2. Graph spaces and universal line bundles.} Let $X$ be
a compact projective manifold. For $d \in H_2(X)$ we call the {\em graph space}
and denote $G_d(X)$
the moduli space of genus $0$ holomorphic maps to $\CC P^1\times X$ of
degree $(1,d)$. This space can be considered as a compactification of 
the space of degree $d$ holomorphic maps $\CC P^1 \to X$. The automorphism
group of $\CC P^1$ acts naturally on $G_d(X)$.
     
Consider first the case $X=\CC P^{N-1}=(\CC ^N-0)/(\CC -0)$. 
Denote $L_d, d\in \ZZ _+$, the projectivization 
$Proj (\CC ^N\otimes S^d (\CC ^2))$ of the
space of degree $d$ vector binary forms 
$(P_1(\z : \xi ): ... : P_N(\z :\xi ))$.
It is a projective space of dimension $Nd +N-1$ and carries a natural
action of $Aut(\CC P^1)=PSL_2(\CC )$ too.

We define the $Aut(\CC P^1)$-equivariant map $\m : G_d(\CC P^{N-1})\to L_d$. 
Consider a stable degree $(1,d)$ map $\psi : C\to \CC P^1\times \CC P^{N-1}$.
There exist a unique irreducible component $C_0\in C$ such that $\psi |_{C_0}$
has degree $(1,d')$ where $d'\leq d$. The image $\psi (C_0)$ is the graph of
a map $\CC P^1\to \CC P^{N-1}$ of degree $d'$. The map is given by the binary
forms $(P_1:...:P_N)$ of degree $d'$ with no common factors and determines the
forms uniquely up to a non-zero constant factor. The curve $C-C_0$ has $r$
connected components which are mapped to $\CC P^1\times \CC P^{N-1}$ with
degrees $(0,d_1),...,(0,d_r), \ d_1+...+d_r=d-d'$, and the image of $i$-th
component is situated in the slice $(a_i :b_i)\times \CC P^{N-1}$.
We put $\m ([\psi ])=\Pi_{i=1}^r (a_i\xi -b_i\z)^{d_i} (P_1:...:P_r)$.    

\medskip

{\bf Proposition 2.1} ( see the Main Lemma in \cite{Gi3}). \newline
{\em The $Aut(\CC P^1)$-equivariant 
map $\m : G_d(\CC P^{N-1})\to L_d$ is regular.}

\medskip

For a compact projective submanifold $X\subset \CC P^{N-1}$ the 
graph space $G_d(X)$ is embedded into the space $G_D(\CC P^{N-1})$
where $D=p(d)$ is the intersection index of the class $p\in H^2(X)$ 
of hyperplane sections with the fundamental class of degree $d$ curves.
Consider the hyperplane line bundle over the projective space $L_D$ and 
induce it to $G_d(X)$ by the restriction of the map $\m $. We will call it
the {\em universal} line bundle corresponding to the embedding of $X$ to 
$\CC P^{N-1}$. The universal line bundle is $Aut(\CC P^1)$-equivariant.

Let now $X$ be a compact projective manifold provided with a 
holomorphic action of a compact Lie group $G$. Let us assume for simplicity
that the $1$-st Chern class identifies the Picard group of $X$ with $H^2(X)$.
Represent a basis in $H^2(X)$ over $\QQ $ by the $1$-st Chern classes of
very ample $G$-equivariant line bundles over $X$.
The embeddings of $X$ to projective spaces defined by holomorphic sections 
of these line bundles are $G$-equivariant. The universal line bundles over
$G_d(X)$ corresponding to these embeddings are $Aut(\CC P^1)\times G$
-equivariant in this case. We can now extend the construction of universal
line bundles from ample to arbitrary elements of the Picard group by 
additivity. For a given equivariant line bundle over $X$ with the 
$G$-equivariant $1$-st Chern class $p\in H^2_G(X)$ we will denote with the same
letter $p$ the $Aut(\CC P^1)\times G$-equivariant $1$-st Chern class
of the corresponding universal line bundle over $G_d(X)$. We will call
$p\in H^2_{PSL_2(\CC )\times G}(G_d(X))$ the {\em universal} class 
corresponding to $p\in H^2_G(X)$.

In the sequel we will use $S^1$-equivariant cohomology with respect to the
maximal torus $S^1\subset PSL_2(\CC )$. We will use the lifting of this 
action to $\CC ^2$ and to the spaces of binary forms defined in the way 
breaking the $PSL_2(\CC )$-symmetry:
\[ (\z , \xi ) \mapsto (\exp (2\pi i \phi )\ \z,\ \xi ) .\]  
If $-p$ denotes the equivariant $1$-st Chern class of the $S^1$-equivariant
Hopf bundle defined by this lifting, then the equivariant integration over
$\CC P^1$ and $L_D$ is described respectively by the residue formulas
\[ f(p,\h )\mapsto \frac{1}{2\pi i} \oint \frac{f(p,\h ) dp}{p (p+\h )} ,\]
\[ f(p,\h )\mapsto \frac{1}{2\pi i} \oint \frac{f(p,\h ) dp}
{p(p-\h )(p-2\h ) ... (p-D\h )} \ ,\]
where $\ZZ [\h ]=H^*(\CC P^{\infty })=H^*(BS^1)$ stands for the coefficient
ring of the $S^1$-equivariant cohomology theory.  

\medskip
{\em Remark.} The graph spaces can be considered as approximations to
the space of loops in $X$. In this interpretation the universal classes
correspond to $S^1$-equivariant forms 
\[ \text{(symplectic structure)}\ +\ \h \ \text{(action functional)} \]
on the loop space. Though heuristicly important, such a relation with 
the Floer theory on the loop space (see \cite{Gi1}) is technically
avoidable. The same applies to Proposition $2.1$ 
(see the remark after Proposition $4.1$).  

\medskip

{\bf 3. Symplectic toric manifolds.} 
Consider the standard real symplectic space $\CC ^N$ provided with
the symplectic structure $\Im \sum d\bar{z}_j\w dz_j/2$. The {\em maximal}
torus $T^N$ acts on $\CC ^N$ by linear symplectic transformations
$z\mapsto diag(\exp 2\pi ix_1 ,..., \exp 2\pi ix_N)z$. The momentum map $J$ 
of this action maps $\CC ^N$ onto the closed $1$-st orthant $\RR _+^N$ in
$Lie ^*T^N$.

Consider a subtorus $T^k\subset T^N$. The momentum map of the action
$T^k : \CC ^N$ is the composition of $J$ with the linear integral projection
$M: \RR^N_+\to \RR ^k=Lie^*T^k$. Pick a regular value $t $ of the momentum
map and define the symplectic orbifold $X$ as the symplectic reduction
$X=\CC //_t T^k=(M \circ J)^{-1}(t )/T^k $. 

Denote $K\subset \RR^k$ the connected component containing $t$ of the regular 
value set. The orbifold $X$ is canonically identified with the complex
quotient $(M \circ J)^{-1}(K)/T^k_{\CC}$. This quotient is a 
Gorenstein projective toric variety of dimension $N-k$. In particular, it has
no codimension $1$ singularities.

The symplectic form on $X$ depends on the choice of the momentum value 
$t\in K$ and is K\"ahler with respect to the complex structure. The action 
of the
quotient torus $T^N/T^k$ on $X$ preserves the form. The momentum polyhedron
of this action is the fiber $M^{-1}(t)$ in the orthant $\RR^N_+$.  The orbifold
$X$ is compact if and only if the polyhedron is compact or, 
equivalently, if $M^{-1}(0)=\{ 0\} $. 
If it is the case, we call $X$ a compact
symplectic toric variety.  

Vice versa, a compact projective toric variety $X$ with at most 
Gorenstein singularities can be obtained as the above symplectic reduction
$\CC^N//_t T^k$ with $N$ equal to the number of hyperplane walls of the
the momentum polyhedron. The polyhedron is described by $N$ inequalities
$H_j\geq c_j $ where $H_j$ are reduced integral linear functions on 
$Lie^*T^{N-k}$. These functions define an embedding of the polyhedron onto 
a section of the $1$-st orthant in $\RR^N$.  

In the sequel we will assume that $X$ is a compact symplectic toric variety
$\CC //_t T^k$ with the number of hyperplane walls of the momentum polyhedron
equal to $N$. We review below some basic properties of $X$ (see \cite{Au, Gi}).

(1) The Gorenstein variety $X$ is non-singular if and only if all the 
$k$-dimensional faces of the orthant $\RR^N_+$ whose projection to $\RR^k$
contain $K$ are mapped to $\RR^k$ with the determinant $\pm 1$.

(2) The correspondence between the regular values $t\in K$ of the momentum map
and cohomology classes of K\"ahler symplectic forms is linear and extends to an
isomorphism $\RR^k\to H^2(X,\RR )$. The isomorphism identifies $K$ with 
the K\"ahler cone of $X$, and the Picard group $H^2(X,\ZZ )$ --- with the
lattice $\ZZ^k\subset \RR^k $ of characters of the torus $T^k$. 
The $1$-st Chern class of the tangent sheaf $\calt _X$ is represented
by the projection $M(u_1+...+u_N)$ of the vector $(1,...,1)\in \RR^N$.

(3) The cohomology algebra $H^*(X,\CC )$ is multiplicatively generated by
K\"ahler classes and can be described as follows. In the space $Lie T^N_{\CC}$
dual to $\RR^N_{\CC }$ consider the union $X^*$ of $(N-k)$-dimensional
coordinate subspaces --- the  orthogonal complements to those 
$k$-dimensional faces of $\RR^N_+$  whose projections to $R^k$ contain $K$.
Denote $\calj$ the ideal of $X^*$ in the algebra $\CC [u_1,...,u_N]$ of 
regular functions on $Lie T^N_{\CC }$. In the algebra $\CC [p_1,...,p_k]$
of regular functions on $Lie T^k_{\CC }$ consider the ideal $\cali $ induced 
from $\calj$ by the embedding $Lie T^k_{\CC }\subset Lie T^N_{\CC }$. Then
$H^*(X,\CC )$ is canonically isomorphic to $\CC [p]/\cali $. The isomorphism
is induced by the correspondence between the infinitesimal characters $u_j$
and $1$-st Chern classes of invertible sheaves on $X$. This description
is valid over $\QQ $, and for non-singular $X$ --- over $\ZZ $.

(4) The algebra $\CC [u]/\calj $ of regular functions on $X^*$ is similarly
identified with the {\em equivariant} cohomology algebra of $X$ with respect
to the action of the quotient torus $T^N/T^k$. By definition, the equivariant 
cohomology algebra of a $G$-space $Y$ is the cohomology algebra of the 
homotopic quotient $EG\times _G Y$ and has the module structure over the 
ring $H^*(BG)$ of $G$-characteristic classes. In the case of a
torus $G$ the ring is identified with the polynomial algebra on $Lie G$. 
The $\CC [Lie T^N/T^k]$-module structure on $\CC [X^*]$ is defined 
by the projection of $X^*$ along $Lie T^k $.

We will use the $T^N$-equivariant cohomology algebra
$H^*_{T^N}(X, \CC )$. Denote $\CC [\l_1,...,\l_N]$ another copy of the 
characteristic class algebra $H^*(BT^N,\CC )$ which will play the role of 
the coefficient ring of $T^N$-equivariant theory throughout the paper. 
Denote $(m_{ij})_{i=1}^k\ _{j=1}^N$ the matrix of the projection
$M:\RR^N_+\to \RR^k$. The $\CC [\l ]$-algebra $H^*_{T^N}(X,\CC )$ is
multiplicatively generated by $(p_1,...,p_k)$ satisfying the relations
$u_j=\sum_{i=1}^k p_i m_{ij} - \l_j$ where $(u_1,...,u_N)$ are the 
generators of $\CC [X^*]=\CC [u]/\calj $. For example, if $M=(1,...,1)$
then $X=\CC P^{N-1}$, $\calj=(u_1...u_N)$, and $H^*_{T^N}(X,\CC )$ is
isomorphic to $\CC [p,\l]/((p-\l_1)...(p-\l_N))$. 
Here $-p$ represents the equivariant $1$-st Chern class of the Hopf bundle
over $\CC P^N$ provided with a natural lift of the torus action. In general
equivariant $1$-st Chern classes of $T^N$-equivariant invertible
sheaves on $X$ are represented by integral linear combinations of $u_j$.
In particular, $u_1+...+u_N$ 
represents the equivariant anti-canonical class of $X$.

(5) Integration over the fundamental cycle  
defines an $H^*(BG)$-linear functional with values in $H^*(BG)$
and the corresponding $H^*(BG)$-bilinear 
Poincare pairing $\lan f, g\ran =\int_X f\w g$ on the equivariant cohomology 
algebra of a compact $G$-manifold. The Borel 
fixed point localization theorem gives rize to the following explicit
description of such integration in the case of the toric symplectic
variety $X$:
\[ \int _X f(p,\l) = \sum _{\a } \Res _{\a} \frac{f(p,\l )\ dp_1\w ...\w dp_k}
{u_1(p,\l)...u_N(p,\l)} \ .\]
Here the index $\a $ runs the set of fixed points of the torus action on $X$.
The fixed points correspond to the vertices of the momentum polyhedron
$M^{-1}(t)$ and therefore --- to those $k$-faces of $\RR^N_+$ whose
projections to $\RR ^k$ contain $t$. The orthogonal complement of such a face
is given by $k$ equations $u_{j_1}=...=u_{j_k}=0, \ j_1<...<j_k$. 
The symbol $\Res _{\a}$ refers to the residue of the $k$-form 
at the pole specified by the ordered set of equations
\[ \sum _{i=1}^k p_i m_{ij_s} = \l_s ,\ s=1,...,k .\]
Permutations of the equations affect the sign of the residue.

We will 
denote $p(\a)=(p_1(\a),...,p_k(\a))$ the solution to this system of linear
equations and $u_j(\a)$ the values of the linear functions
$u_j=\sum p_im_{ij}-\l_j$ at the point $p(\a)$.
We will also identify in our notations the index $\a $ of the fixed point
of $T^N : X$ with the corresponding vertex of the momentum polyhedron and with
the set $\{ j_1, ... , j_k \} $.

Although the residues on the RHS are rational functions of $\l $, their sum
is a polynomial provided that $f\in H^*_{T^N}(X,\CC )$ is represented by
a polynomial $f(p,\l)$. The equivariant Poincare pairing 
$\lan f, g \ran = \int _X fg $ is
non-degenerate over $\QQ [\l ]$ (and even $\ZZ [\l]$ if $X$ is non-singular). 
The value at $\l=0$ of the LHS represents $\int _X f(p,0)$ and gives rise 
to the ordinary Poincare pairing on $H^*(X)$.    

Notice that the toric variety $X$ and its properties have been
described entirely in terms of the $N$ integer vectors $Mu_1,...,Mu_N$
in $\RR ^k$ and the chosen K\"ahler cone $K\subset \RR^k$.

\medskip

We will assume further on that $X$ denotes a {\em non-singular} compact 
symplectic toric manifold.

\medskip

A codimension $l$ complete intersection $Y\subset X$ is by definition the 
common zero locus of holomorphic sections of $l$ line bundles. 
Let $v_1,...,v_l$
be the infinitesimal characters of $T^N$ which
specify $l$ equivariant line bundles $\call_a$. 
We will assume that $\call_a$ are
non-negative, \this $M v_a\in \bar{K}$, and therefore that the bundle 
$\calv =\oplus _a \call_a $ is convex.

We provide the bundle $\calv $ with the additional action of the 
$l$-dimensional 
torus $T^l$ acting fiberwise by scalar multiplication in each summand 
$\call_a $. We will study the equivariant Gromov - Witten 
theory of the pair ($X$, $\calv $) with respect to the action of the torus
$T=T^N\times T^l$. The algebra
\[ \QQ [ \l_1,...,\l_N, \l'_1,...,\l'_l ]=H^*(BT,\QQ ) \] 
will play the role of the 
coefficient algebra of the theory. Now on we will assume that $u_j, v_a, p_i$
denote corresponding $T$-equivariant cohomology classes. In particular the
$T$-equivariant $1$-st Chern classes of the line bundles $\call_a$ and their
fixed point localizations are represented by the linear combinations 
$v_a:= \sum l_{ia}p_i-\l'_a$, and $v_a(\a)= \sum l_{ia}p_i(\a)-\l'_a$.

We introduce the $T$-equivariant integration
over the {\em virtual} fundamental class 
$[Y]$ of the corresponding (invariant) complete intersections:
\[ \int _{[Y]} f \ :=\ \int_X f \ Euler (\calv ) = 
\sum _{\a} \Res _{\a} f (p, \l, \l') \ \frac{ v_1 ... v_l \ dp_1\w ... \w dp_k}
{u_1 ... u_N} .\]
In the non-equivariant limit $\l=0, \l'=0$ it reduces to the integration 
over the fundamental class of non-singular comlete intersections $Y\subset X$
as well as the whole Gromov-Witten theory --- to that for $Y$, according to
the axiom (7).
 
\medskip

Degrees of compact holomorphic curves in $X$ form the semigroup 
\[ \L = \{ d\in H_2(X)| (t, d)\geq 0 \ \forall t\in  \bar{K} \} .\] 
The Kahler cone $\bar{K}$ is the intersection $\cap_{\a}\D_{\a} $ of 
of images in $\RR ^k$ of $k$-faces in $\RR_+^N$. Respectively, $\L $ coincides 
with the convex hull of the union $\cup _{\a} \D_{\a}^*$ of the orthants
$\D^*_{\a}\subset \RR^{k*}$ polar to the orthants $\D_{\a}$.
We will denote $\ZZ [[\L ]]$ the formal completion of the semigroup ring 
$\ZZ [\L ]$ and represent an element $d\in \L$ by the monomial 
$q^d=q_1^{d_1}...q_k^{d_k} \in \ZZ [\L ]$ where $(d_1,...,d_k)$ are 
coordinates of $d$ in $(\ZZ^k)^*$.

We will write $d\geq d'$ if $d-d'\in \L$. It is a partial order on $H_2(X)$.

We will also use the notations 
\[ D_j=\sum_i d_i m_{ij}, \ L_a=\sum_i l_{ia}d_i \]
(which mask the actual dependence of these integers on the vector $d$).

\medskip

We complete this section with some notations and elementary information
about the $1$-dimensional orbits of $T_{\CC}^N:X$ 
which will be exploited, in accordance with the axiom (8),
in the fixed point localization technique. 

The $1$-dimensional orbits correspond to the $1$-dimensional edges of the
momentum polyhedron. 
The vertex $\a $ of the momentum 
polyhedron $M^{-1}(t)$ is situated in the $k$-face 
$\{ (T_1,...,T_N)| T_s=0, \forall s\notin \a \} $
of $\RR ^N_+$ and connected by one-dimensional edges to $N-k$ other vertices 
$\b (\a , j), j\notin \a$, situated in the $(k+1)$-faces
$\{ (T_1,...,T_N) |  T_s=0 \forall s\notin \a , s\neq j \} $. The edge 
itself is the momentum polyhedron of the closure of a $1$-dimensional toric 
orbit 
$\CC ^{k+1} //_t T^k $ in $X$ isomorphic to $\CC P^1$. We denote 
$d(\a , j)\in \L $ the degree of this $\CC P^1$ in $X$. The corresponding
$D_s(\a , j):=\sum d_i(\a,j)m_{is}$ vanish for $k-1$ values of the index 
$s\in \a \cap \b (\a,j)$ and are equal to $1$ for the two values of 
$s\in (\a \triangle \b)=\{ j, j'\} $.  We denote 
$L_a(\a , j)$ the corresponding values $\sum d_i(\a ,j) l_{ia}$ of $L_a$.   
Notice that $\b $ and $j'\in \a $ depend on and are uniquely determined by a
choice of $\a $ and $j\notin \a $.

For $\b =\b(\a,j)$ we have $d(\a ,j)=d(\b ,j')$ and
\[ u_s(\a)=u_s(\b)+D_s(\a,j)u_j(\a), \ 
v_a(\a) = v_a(\b ) +L_a(\a, j)u_j(\a), \]
In particular, $u_j(\a)=-u_{j'}(\b)$. 

\medskip

{\bf 4. The generating function.} 
Let $(X,\calv )$ be the nonsingular compact toric symplectic manifold provided
with the convex bundle equivariant with reapect to the action of 
$T=T^N\times T^l$ as described in the previous section.
Denote $\calv _d$ the vector bundle over $G_d(X)$ with the fiber
$H^0(C,\psi ^* \pi_2 ^* \calv )$ over the point represented by the stable map
$\phi :C\to \CC P^1\times X$.  We will study the following series:
\[ \calg (z,q,\h; \l, \l' ) =\sum _{d\in \L} q^d 
\int _{[G_d(X)]} e^{P_1z_1+...+P_kz_k} Euler (\calv_d) .\]
Here $(P_1,...,P_k)$ denote the universal classes in 
$H^2_{T\times S^1}(G_d(X),\QQ )$ corresponding to the $T$-equivariant
classes $p_i\in H^2_T(X)$.

\medskip

{\bf Proposition 4.1.} {\em The series $\calg = \sum  g_{d,m} q^d z^m$ has 
polynomial coefficients $g_{d,m}\in \QQ [\h ,\l,\l' ]$.}

{\em Proof.} The equivariant cohomology classes $P_i$ and 
$Euler (\calv_d)$
of $G_d(X)$ are defined over $\QQ [\h ,\l,\l' ]$, and the integration over 
the Gromov-Witten class $[G_d(X)]$ assumes values in $\QQ [\h ,\l,\l' ]$ by the
axiom (1). 

\medskip

{\em Remark.} The polynomiality property of $g_{d,m}$ is crucial for our
proof of Theorems $0.1$, $0.2$. However one can define the series $\calg $
in terms of Gromov-Witten theory on $\CC P^1\times X$ (see \cite{Gi3} , 
Section $6$) without mentioning the universal classes $P_i$, and the 
polynomial property then follows directly from the axiom $(1)$. Thus our 
use of Proposition $2.1$ is avoidable.
   
\medskip

{\bf Proposition 4.2.} (see \cite{Gi3}). 
\[ \calg ( z, q, \h )=
\int_{[Y]}  \cals (q \exp (\h z), \h)\ e^{pz}\ 
\cals (q , -\h)\ , \]
{\em where $\cals (q, \h; \l, \l')$ is determined from the condition that 
for any $a\in H^*_T(X)$}
\[ \int_{[X]} \cals \ a =
\int_{[X]} Euler (\calv )\ a  + \sum _{d\in \L-0} q^d \int_{[X_{2,d}]}
\frac{\ev_1^* (a) \ Euler (\calv_{2,d})}{(\h - c)} \ .\]

{\em Proof.} The proposition is easily deduced (see \cite{Gi3}) by 
localization to fixed points of the $S^1$-action
on $G_d(X)$ and $L_D$ and from the following corollary of the axiom (2):
\[ \int _{[X_{2,d}]} \frac{\fgt_2^*(A)}{\h - c} = \int _{[X_{1,d}]} 
\frac{A}{\h (\h - c)} . \] 

{\bf Corollary 4.3.} {\em $\cals = 1 + o (1/\h ) $ when written as
a formal power series in $\h ^{-1}$.}

{\em Proof:} $\ Euler (\calv_{2,d})=\fgt_2^* (Euler (\calv_{1,d}))$.

\medskip

{\em Remark.} The vector-function $J_{\calv }$ in Theorem $0.2$ 
is defined in Section $1$ as  $\cals \ \exp ((t_0+\sum P_i\log q_i)/\h ) $.

\medskip

Denote $\cals_{\a }$ the restriction of the $T$-equivariant cohomology class 
$\cals =\sum _{d\in \L} S_{\a, d} q^d $
to the fixed point $\a $ of the torus $T$ action on $X$. The coefficients
$\cals _{\a, d}$ of the series $\cals _{\a}$ are rational functions of 
$(\h ,\l,\l')$.

\medskip

{\bf Proposition 4.4.} {\em The series $\cals_{\a}$ satisfy the following 
recursion relations:
\[ \cals_{\a }(q,\h;\l,\l')=\sum_{d\in \L} R_{\a, d}(\h^{-1};\l,\l' ) q^d + \]
\[ \sum_{j\notin \a } \sum_{n=1}^{\infty } q^{n\ d_{\a,j}}
\frac{C_{\a, j}(n)}{n\h+u_j(\a )} \cals _{\b (\a, j)} (q, -u_j(\a)/n; \l,\l' ) 
\ ,\]
where $R_{\a ,d}$ are polynomials in $\h^{-1}$ with coefficients in
$\QQ (\l,\l' )$, and the recursion coefficients $C_{\a ,j}(n)=$} 
\[ \frac{\Pi_a\Pi_{m=1}^{nL_a(\a,j)} (v_a(\a )-mu_j(\a )/n)
\ \Pi_{s\notin \a \cup \b }\Pi_{m\leq 0}(u_s(\a)-mu_j(\a)/n))}
{(n-1)! (u_j(\a)/n)^{n-1})\ 
\Pi_{s\notin \b }\Pi_{m\leq nD_s(\a ,j)} (u_s(\a )-mu_j(\a)/n)}\ . \]

{\em Proof.} We apply the fixed point localization technique to the action
of the torus $T$ on $X_{2,d}$. Consider a fixed point represented by the 
stable map $f:(C, x_1, x_2) \to X$ and denote $C'$ the irreducible component of
$C$ carrying the marked point $x_1$. The contribution of this fixed point to
${\cal S}_{\a}$ via the localization formula described in the axiom $(8)$ 
vanishes unless $f(x_1)=\a $. 

If $f(C')=\a $, then the contribution involves
integration over the space $\bar{\frak M}_r$ of configurations of $r>2$ special
points on the connected component of $f^{-1}(\a)$ containing $x_1$. The
universal cotangent line at the $1$-st marked point localizes to the 
line bundle over $\bar{\frak M}_r$ formed by the cotangent lines to $C'$ at 
$x_1$. Since the action of $T$ on this line bundle is trivial, the 
localization of the Chern class $c$ is nilpotent. Thus the contribution of such
fixed points to the localization formula is polynomial in $1/\hbar $.

If $f(x_1)=\a$ but $f(C')\neq \a $ then $f':=(f | C')$ is the degree $n$ 
cover of a $1$-dimensional orbit of $T_{\CC }$ in $X$ 
connecting the fixed point $\a $ with another fixed point $\b =\b (\a,j) $. 
In this case $f$ is glued from $f':(C',x_1,x)\to X$ and a stable map 
$f'':(C'',x,x_2)\to X$ of degree $d''=d-nd(\a, j)$ with $f'(x)=f''(x)=\b $.
The localization of $\h -c$ in this case equals $\h + u_j(\a)/n$. The 
contribution of $T'_x\otimes T''_x$ to the Virtual Normal Euler Class in the 
localization axiom equals $-u_j(\a)/n -c$ where $c$ is the equivariant $1$-st
Chern class of the universal cotangent line over $X_{2, d''}$ at the marked
point $x$. The remaining part of the Virtual Normal Euler Class is the product
of such classes for $f''$ in $X_{2,d''}$ and for $f'$ in $X_{2, nd(\a ,j)}$.
The latter one can be easily computed by using the quotient description
of the tangent bundle to $X=\CC^N//T^k $ as a virtual direct sum of line 
bundles and is equal to $1/C_{\a, j}(n)$. Also we have
$|Aut f|= n |Aut f''|$ where $n$ is the order of the cyclic automorphism
group of the $n$-fold cover $f': z\mapsto z^n$.

The summation over all $d''$, $ j\notin \a$ and $n=1,2, ... $ gives rise to
the recursion relation. $\square $ 

\medskip

{\bf Proposition 4.5.} {\em The series $\cals $ is uniquely determined by the
following properties:

(a) the recursion relations of Proposition $4.4$,

(b) the asymptotical condition $\cals = 1+o (1/\h ) $ of Corollary $4.3$,
 
(c) the property from Proposition $4.1$ of the series 
\[ \sum g_{d,m}q^dz^m =
\int_{[Y]} \cals (qe^{\h z},\h) e^{pz} \cals (q,-\h) \]
to have polynomial coefficients $g_{d,m}\in \QQ [\h,\l,\l' ]$.}

{\em Proof.} Consider another solution $\cals '$ to the same recursion
relations (a) with the same polynomials $R'_{\a, d}=R_{\a, d}$ as in $\cals $ 
for all $0\leq d<d_0$ and satisfying the asymptotical condition (b). 
We will show that $R'_{\a,d_0}=R_{\a, d_0}$. 
This would imply the uniqueness since, given
all the polynomials $R_{a,d}$, the recursion relations (a) allow one to
recover the series $\cals $ unambiguously.

Consider the $q^{d_0}$-term in $\calg ' -\calg $. The conditions 
(a),(b) imply that $\cals '-\cals = R q^{d_0} +\ 
(\text{higher order terms in} \ q)$
and that $q^{d_0}$-term in $\calg '- \calg$ is equal to
\[ \d (R)=\int_{[Y]} e^{(p+d_0\h)z} R(\h,\l ) + e^{pz} R(-\h,\l ) \ ,\]
where the class $R$ is defined by its localizations 
$R_{\a}=R'_{\a ,d_0}-R_{\a ,d_0}$. From (b) and (c) we know that $R$ is
a polynomial in $1/\h $  divisible by $1/\h ^2$ and that
the series $\d =\sum \d_m z^m $ has polynomial coefficients 
$\d _m\in \QQ [\h ,\l,\l' ]$.

A generic value of $\l =(\l_1,...,\l_N)$ determines a linear function
on the momentum polyhedron $\RR _+^N\cap J^{-1} (z)$ for each $z\in K$
with pairwise distinct values $H_{\a}$ at the vertices. Computing 
$\d ( A\h^{-1}+B )$ modulo $\h $ for a generic ray $t\mapsto tz $ 
we find
\[ \sum_{\a} e^{H_{\a} t} (A_{\a}t\ (z,d_0)  +B_{\a} )
\Res _{\a} \frac{v_1...v_l \ dp_1\w ... \w dp_k}{u_1...u_N} .\]
Here $v_1(\a)...v_a(\a)\neq 0$ for generic values of $\l'$, and $(z,d_0)>0$.
Since the functions $\exp (H_{\a} t), t\exp (H_{\a}t)$ with distinct $H{\a}$
are linearly independent, we conclude that 
$\d (A/\h +B ) =o (\h )$ implies $A=B=0$ in $H^*_T(X, \CC (\l,\l'))$
for generic and therefore for all values of $\l, \l' $ and $z$.
Applying this conclusion to $\h^{2r} R(1/\h, \l )$ with $r>0$ and
$R(1/\h )= A\h^{-2r-1}+B\h^{-2r}+...$ of degree $\leq (2r+1)$ we find that
$\deg R\leq (2r-1)$ and by induction --- that $\deg R \leq 1$. Now the
assumption (b) that $R$ has no terms of order $\h ^{0}$ and $\h ^{-1}$
implies that $R=0$. $\square $
 
\medskip

{\bf 5. Toric map spaces.} 
In this section we describe toric compactifications of spaces of 
holomorphic maps $\CC P^1 \to X$. 

For each $d\in \L $ consider the following complex 
space of $N$-dimensional vector-polynomials in one complex variable $\z$:
\[ \{ (z_1(\z),...,z_N(\z)) | \deg z_j\leq D_j \} \]
The torus $T^N$ acts on this space componentwise. Denote $J_d$ the momentum map
of the induced action of $T^k\subset T^N$, pick $t\in K$ and introduce the 
symplectic toric variety $X_d=J_d^{-1}(t)/T^k$. It is compact and nonsingular 
whenever $X$ is compact and nonsingular. Generic points in $X_d$ represent
degree $d$ holomorphic maps 
\[ \CC P^1\to X: \z \mapsto (z_1(\z),...,z_N(\z)) \mod  T^k_{\CC }\ .\]
The variety $X_d$ is empty unless $d\in \cup_{\a}\D^*_{\a}$.

The rotation $\z\mapsto \exp (2\pi i \phi) \z$ of $\CC P^1$ induces an
$S^1$-action on $X_d$. Fixed points of the $S^1\times T^N$-action on $X_d$
are isolated. The Borel localization
formula for $X_d$ yields:
\[ \int _{X_d} f(p,\h,\l) = \sum _{\a: d\in \D^*_{\a}} \sum_r
\Res _{\a, r} \frac{f(p,\h,\l) dp_1\w ...\w dp_k}
{\Pi _{j: D_j\geq 0} (u_j(u_j-\h)(u_j-2\h)...(u_j-D_j\h))}\ .\]
Here $\Res _{\a ,r }$ refers to the residue at the point specified by the
equations 
\[ u_{j_1}(p,\l )=r_1\h \ ,...,\ u_{j_k}(p,\l )=r_k\h \ , j_1<...<j_k, \]
with $\{ j_1,...,j_k\} = \a $ and $r$ runs the integer vectors 
\[ r=(r_1,...,r_k) : 0\leq r_1 \leq D_{j_1}, ..., 0\leq r_k\leq D_{j_k} \ .\]
The equivariant cohomology algebra of $X_d$ is identified with the quotient
of $\CC [p,\h,\l ]$ by the kernel ideal of the corresponding Poincare pairing.

The dimension of $X_d$ may exceed the Riemann-Roch dimension 
$\dim X+ (c_1(\calt_X), d)=N-k+D_1+...+D_N$ 
of the space of degree $d$ holomorphic maps $\CC P^1\to X$ if some 
$D_j=\sum d_i m_{ij} $ are
negative. We introduce the {\em virtual fundamental class} $[X_d]$ of 
the Riemann-Roch dimension as the Poincare-dual to the equivariant Euler 
class of the following vector bundle over $X_d$.

Consider the complex space of the vector polynomials in the variable $\z^{-1}$
defined by
\[ \{ (z_1(\z^{-1}),...,z_N(\z^{-1}))| z_j(0)=0, \deg z_j< -D_j\} \ .\]
We introduce the $T^N$-equivariant locally free sheaf on 
$X_d=J_d^{-1}(t)/T^k$ associated with the componentwise action of 
$T^k\subset T^N$ on this space. The rotation of $\z $ defines an
$S^1$-equivariant structure on the sheaf. 
At generic points $\psi: \CC P^1\to X$ of $X_d$ the sheaf coincides with the
obstruction sheaf $H^1(\CC P^1, \psi ^*(\calt _X))$.

The Borel localization formula gives rise to the following description of 
integration over the virtual fundamental class:
\[ \int _{[X_d]} f = \sum _{\a: d\in \D^*_{\a}} \sum_r
\Res_{\a , r} dp_1\w ... \w dp_k \ f(p,\h,\l ) 
\Pi_{j=1}^N \frac{\Pi _{m=-\infty }^{-1} (u_j-m\h )}
            {\Pi _{m=-\infty }^{D_j} (u_j-m\h )} \ .\]    

\medskip

Consider now the convex bundle $\calv = \oplus_a \call_a $.
For a generic degree-$d$ map $\psi :\CC P^1\to X$ the space 
$H^0(\CC P^1, \psi ^* \call_a)$ can be identified with the space of polynomials
$ \{ y(\z) | \deg y \leq L_a \} $ where $L_a=\sum d_i l_{ia} $. 
We introduce the vector bundle of dimension $L_a+1$ over $X_d$ associated
with the scalar action of $T^k$ on this space via the character 
$\sum l_{ia}p_i$. The direct sum of these bundles over $a=1,...,l$ is 
equivariant with respect to $S^1\times T^N\times T^l$ where $S^1$ acts by
rotations of $\z $, and $T^l$ acts by scalr multiplication 
componentwise on the direct summands. The $S^1\times T$-equivariant 
Euler class of this bundle is 
$\Pi _a\ v_a (v_a -\h )...(v_a -L_a\h )$ where $v_a=\sum l_{ia}p_i-\l'_a$. 
We define the 
{\em virtual fundamental class} $[Y_d]$ as the Poincare-dual to 
this Euler class:
\[ \int _{[Y_d]} f := \int _{[X_d]} f(p,\h, \l, \l' ) 
\Pi _a \Pi _{m=0}^{L_a} (v_a(p,\l )-m\h ) \ .\]

For polynomial $f\in \QQ [p,\h,\l,\l']$ the integration over $[Y_d]$ 
assumes polinomial values in $\QQ [\h,\l,\l']$.

\medskip

{\bf 6. The hypergeometric series.} 
We will study here the generating function
\[ \Phi (z,q) = \sum_{d\in \L} q^d 
\int _{[Y_d]} e^{p_1z_1+...+p_kz_k} \]
which is a formal power series in $z$ with coefficients in the formal
completion of the semigroup ring of $\L $ over the polynomial algebra 
$\QQ [\h,\l,\l']$. 

\medskip

{\bf Proposition 6.1.} {\em The series $\Phi $ is weighted-homogeneous of 
degree $l+k-N$ with respect to the grading}
\[ \deg z_i=-1, \deg \h =\deg \l_j=\deg \l'_a =1,
 \deg q_i = \sum_j m_{ij} - \sum_a l_{ia} .\]

{\em Proof:} 
This follows from grading in equivariant cohomology and the formula
$N-k-l + (c_1(\calt _X), d)- (c_1(\calv ),d)$ for the complex dimension
of the virtual fundamental cycle $[Y_d]$.

\medskip

Consider the formal $q$-series $\Psi $ (which differs from $I_{\calv}$ in
Theorem $0.2$ by the exponential factor $\exp ((t_0+p\log q)/\h )$ only):
\[ \Psi (q,\h;\l,\l' )=\sum_{d\in \L } q^d 
\frac{\Pi_a \Pi_{m\leq L_a} (v_a+m\h ) \Pi_j \Pi_{m\leq 0} (u_j+m\h)}
{\Pi_a \Pi_{m\leq 0} (v_a+m\h ) \Pi_j \Pi_{m\leq D_j} (u_j+m\h )} \ ,\]
where $u_j=\sum p_i m_{ij} -\l_j$ and $v_a=\sum p_i l_{ia}-\l'_a$ are 
$T$-equivariant 
cohomology classes of $X$, the integers $D_j=\sum d_i m_{ij}$,
$L_a=\sum d_i l_{ia}$ depend on $d=(d_1,...,d_k)$, and $m$ runs integer
values starting from $-\infty $. 
Coefficients of the series $\Psi$ 
are well-defined equivariant cohomology classes
of $X$ over the field $\CC (\h ,\l )$ of rational functions, and the whole 
series can be considered as an equivariant cohomology class of $X$ with
appropriate coefficients. Notice that a $q^d$-term in $\Psi $ has zero
localization at the fixed point $\a$ unless $d\in \D^*_{\a}$.

\medskip

{\bf Proposition 6.2} (compare with \cite{Gi1,Gi2,Gi3}).
\[ \int _{[Y]} \Psi (q\exp (\h z) ,\h;\l,\l' ) e^{pz} \Psi (q,-\h; \l,\l' ) =
\Phi (z, q,\h ,\l,\l') .\]
{\em In particular, $\Phi ((t-\t )/\h, e^{\t },\h )$ is invariant under the
change $(t,\t ,\h )\mapsto (\t ,t,-\h )$.}

{\em Proof.} We have
\[ \int_{[Y]} \Psi (q\exp (\h z),\h) e^{pz} \Psi (q,-\h) =\]
\[ \sum_{\a} \sum_{d',d''\in \D^*_{\a}} q^{d'+d''} \Res_{\a} 
e^{(p+\h d')z}\  dp_1\w ...\w dp_k\ \times \]
\[ \frac{\Pi_a \Pi_{m=-L''_a}^{L'_a} (v_a+m\h)\ 
\Pi_j [\Pi_{m\leq 0}(u_j+m\h) \ \Pi_{m\geq 0} (u_j+m\h) ]}
{ \Pi_j u_j [\Pi_{m\leq D'_j} (u_j+m\h)\ \Pi_{m\geq -D''_j} (u_j+m\h)]} \]
\[ =\sum_{\a} \sum_{d\in \D^*_{\a}} q^d \sum_r \Res_{\a,r} 
 e^{pz} dp_1\w ...\w dp_k\ \times \]
\[ \frac{\Pi_a \Pi_{m=-L_a}^0(v_a+m\h)\ 
\Pi_j u_j \Pi_{m\in \ZZ} (u_j+m\h)}
{\Pi_j u_j [\Pi_{m\leq 0} (u_j+m\h) \ \Pi_{m\geq -D_j} (u_j+m\h)} \]
\[ = \sum_{\a} \sum_{d\in \D^*_{\a}} \sum_r \Res_{\a,r}
\frac{e^{pz}dp_1\w ...\w dp_k\ \Pi_a \Pi_{m=0}^{L_a} (v_a-m\h)\ 
\Pi_j \Pi_{m<0} (u_j-m\h)}{\Pi_j \Pi_{m\leq D_j} (u_j-m\h)}\]
\[= \sum_{d\in \L} q^d \int_{[Y_d]} e^{pz} .\]
$\square $

\medskip

Consider the localizations $\Psi _{\a} (q,\h ;\l,\l')$ of the class $\Psi$ 
at the fixed points:
\[ \Psi _{\a} = \sum _{d\in \D^*_{\a}} q^d 
\frac{\Pi _a \Pi _{m=1}^{L_a} (v_a(\a)+m\h)}{\Pi _{j\in \a} D_j! \h^{D_j}}\ 
\Pi _{j\notin \a} \frac{\Pi _{m\leq 0} (u_j(\a)+m\h)}{\Pi_{m\leq D_j} 
(u_j(\a)+m\h )} .\]
Each $q^d$-term in this series is a rational function of $\h$ and can be
uniquely written as the sum of simple fractions with poles
at $\h =-u_j(\a)/n, \ j\notin \a,\ 0< n\leq D_j$ plus a Laurent polynomial 
in $\h$.

\medskip

{\bf Proposition 6.3.} {\em The series $\Psi _{\a}$ satisfy the
following recursion relations:
\[ \Psi _{\a} (q,\h; \l,\l')=\sum_{d\in \D^*_{\a}} q^d R_{\a,d} (\h,\l,\l')+\]
\[
\sum_{j\notin \a } \sum_{n>0} q^{n\ d(\a,j)} \frac{C_{\a,j}(n)}{(u_j(\a)+n\h)}
\Psi _{\b } (q, -u_j(\a)/n;\l ,\l' ) \ ,\]
where $\b=\b(\a,j)$ and the coefficients} $C_{\a,j}(n)=$
\[ \frac{\Pi_a \Pi_{m=1}^{nL_a(\a,j)}(v_a(\a)-mu_j(\a)/n) \ 
\Pi_{s\notin \a \cup \b} \Pi_{m\leq 0} (u_s(\a)-mu_j(\a)/n)}
{(n-1)! (u_j(\a)/n)^{n-1}\ 
\Pi_{s\notin \b} \Pi_{m\leq nD_s(\a,j)}(u_s(\a)-mu_j(\a)/n)}\ .\]

{\em Proof.} For generic $\l,\l'$ and $\b =\b (\a,j)$ denote $Z_{\b}^d$ the 
coefficient of the series $\Psi _{\b}$ at $q^d$. 
The value of this coefficient at $\h = -u_j(\a)/n=u_{j'}(\b)/n$ is equal to
\[ \frac{\Pi_a\Pi_{m=1}{L_a}(v_a(\b)-mu_j(\a)/n) \ 
\Pi_s\Pi_{m\leq 0} (u_s(\b)-mu_j(\a)/n)}
{\Pi_s \Pi_{m\leq D_s} (u_s(\b)-mu_j(\a)/n)}\]
It vanishes unless $d\in \D^*(\b)$ and $D_{j'}+n\geq 0$. 
Due to the relations 
$u_s(\b)=u_s(\a)-D_s(\a,j)u_j(\a),\ v_a(\b)=v_a(\a)-L_a(\a,j)u_j(\a)$
we have $Z_{\b}^d(-u_j(\a)/n )=$
\[  \frac{\Pi_a\Pi_{m>nL_a(\a,j)}^{L_a+nL_a(\a,j)} (v_a(\a)-mu_j(\a)/n)\ 
\Pi_s \Pi_{m\leq nD_s(\a,j)} (u_s(\a)-mu_j(\a)/n)}
{\Pi_s \Pi_{m\leq D_s+nD_s(\a,j)} (u_s(\a)-mu_j(\a)/n)} \ .\]
Since $D_s(\a,j)=0$ for $s\in \a \cap \b$ and $D_j(\a,j)=D_{j'}(\a,j)=1$,
we have $C_{\a,j}(n) Z_{\b}^d(-u_j(\a)/n)=$
\[ \frac{\Pi_a\Pi_{m=1}^{L_a+nL_a(\a,j)}(v_a(\a)-\frac{m}{n}u_j(\a))\ 
\Pi_{s\neq j,j'} \Pi_{m\leq 0} (u_s(\a)-\frac{m}{n}u_j(\a))}
{\Pi_{m=1}^{D_{j'}+n} (\frac{-n-m}{n}u_j(\a))\ 
\Pi_{m=1}^{D_j+n}\ _{m\neq n} (\frac{n-m}{n}u_j(\a))\ 
\Pi_{s\neq j,j'} \Pi_{m\leq D_s+nD_s(\a,j)} (u_s(\a)-\frac{m}{n}u_j(\a))} \ .\]
The last product is exactly the residue of the $q^{d+nd(\a,j)}$-term in
the series $\Psi _{\a}$ at the simple fraction with the pole $\h= -u_j(\a)/n$.
The vector $d'=d+n d(\a,j)$ is automatically in $\D^*_{\a}$, and additionally
$D_j+n\geq n$. 
Vice versa, when $d'$ is in $\D^*_{\a}$, and for some $j\notin \a$ and 
$n\in \NN $ we have $D'_j\geq n$, then $d=d'-nd(\a,j)$ is in $\D^*_{\b(\a,j)}$
and $D_{j'}+n=D'_{j'}\geq 0$.
 
Thus the sum of simple fractions on the RHS of the recursion relations
reproduces all the simple fractions of the LHS with non-zero poles, and
the remaining part of each $Z_{\a}^{d'}(\h)$ is a Laurent polynomial of $\h$.
Since the identity between the $q$-series with coefficients in 
$\QQ (\h ,\l, \l')$ holds for generic values of $(\l,\l')$ it holds in fact
over $\QQ (\h, \l, \l')$. $\ \square $

\medskip

{\em Remark.} The coefficients $C_{\a ,j}(n)$ in the recurrence relations of
Propositions $4.4$ and $6.3$ are the same. Although we obtained this fact by 
straightforward computations of both coefficients, the coincidence should
not be considered a miracle. The value of the coefficient 
$C_{\a,j}(n)$ in Proposition $4.4$ is determined by localizations to fixed
points represented by {\em irreducible} degree $n$ maps
$\CC P^1 \to X$. On the other hand, the graph spaces differ from 
the toric map spaces only near stable maps of reducible curves. Perhaps the
coincidence of the recurrence coefficients can be derived from this
observation.

\medskip
  
Consider now the case of complete intersections $Y$ with non-negative
$1$-st Chern class $c_1(\calt_Y)$. This means that
\[ \sum_i p_i (\sum_j m_{ij} - \sum_a l_{ia}) \ \in \ \bar{K} \ .\]
A reformulation of this condition 
reads: $\deg q^d \geq 0$ for all $d\in \L$. 

Put 
\[ \L_0=\{ d\in \L | \sum L_a=\sum D_j, D_j\geq 0 \forall j=1,...,N\} ,\]
\[ \L_1=\{ d\in \L | \sum D_j-\sum L_a =1, D_j\geq 0 \forall j=1,...,N \} ,\]
\[ \L'_0=\{ d\in \L | \sum L_a=\sum D_j \} .\]

\medskip

{\bf Proposition 6.4.} {\em Suppose that $\calt_Y$ is non-negative. Then}
\[ \Psi = \Psi ^{(0)}+\Psi ^{(1)}/\h + o(1/\h ) , \]
\[ \Psi ^{(0)}=\sum_{d\in \L_0}
 \frac{L_1!...L_l!}{D_1!...D_N!} q^d \ ,\]
\[ \Psi ^{(1)}= \sum_{d\in \L _1}\frac{L_1!...L_l!}{D_1!...D_N!} q^d +
\sum _a g_a(q) v_a - \sum_j f_j(q) u_j  ,\]
{\em where $g_a$ and $f_j$ are some power series $\sum_{d\in \L'_0} A_d q^d$
with coefficients $A_d\in \QQ $ and $A_0=0$.}


{\em Proof.} In the case of non-negative $\calt_Y$ we have
$\sum D_j \geq \sum L_a$ for all $d\in \L$. 
The definition of $\Psi (q,\h;\l,\l' )$ then shows that 
it can be written as a power series in $1/\h $, and that 
the terms of order $\h^0$ and $\h^{-1}$ have the form described in the
proposition. 


\medskip

{\bf 7. Equivalence transformations.} 
We study further the case of complete intersections $Y$ with non-negative
tangent sheaf.

Consider a series $Z(q,\h ; \l, \l')\in H^*_T(Y,\QQ)$ of weighted-homogeneous 
degree $0$ with respect to the grading $\deg \h=\deg \l =\deg \l' =1,\ 
\deg q_i=\sum_j m_{ij}-\sum_a l_{ia}$ which satisfies

(i) the recursion relation  
\[ Z_{\a}(q,\h;\l,\l')=\sum_{d\in \L} R_{\a,d}(\h^{-1};\l,\l') q^d +\]
\[ \sum_{j\notin \a} q^{n\ d_{\a,j}} \frac{C_{\a,j}(n)}{n\h+u_j(\a)}
Z_{\b(\a,j)} (q,-u_j(\a)/n;\l,\l')\ , \]
where $R_{\a,d}$ are polynomials in $\h^{-1}$ with coefficients in 
$\QQ(\l,\l')$ and the recursion coefficients $C_{\a,j}(n)$ are described
in Propositions $4.4$ and $6.3$;

(ii) the condition that the series
\[W= \sum w_{d,m} q^d z^m := \int_{[Y]} Z(qe^{\h z},\h )e^{pz}Z(q,-\h) \]
has polynomial coefficients $w_{d,m}\in \QQ[\h,\l,\l']$.

In Section $4$ we proved that the series $\cals $ is uniquely determined by 
these properties and the asymptotical condition of Corollary $4.3$. In the
previous section we found a series $\Psi $ which also has these properties
but may differ from $\cals $ by the asymptotical condition. In this section 
we describe some transformations of the series $Z$ which preserve the 
properties (i),(ii) but change the asymptotical expansion 
$Z=Z^{(0)}+Z^{(1)}/\h +o (1/\h ) $.

\medskip

{\bf Proposition 7.1.} {\em Let $f=\sum_{\L'_0} f_dq^d$ be a series with 
rational coefficients $f_d\in \QQ $ and $f_0\neq 0$. Then $f Z$ satisfies the
conditions (i),(ii).}

{\em Proof.} Simultaneous multiplication of the localizations $Z_{\a}$ by $f$
obviously preserves the the recursion relation of Proposition $4.4$ and 
affects the series $\sum R_{\a,d} q^d$ in such a way that the polynomiality
property of $R_{\a,d}(1/\h)$ is preserved. This proves (i).

The multiplication of $Z$ by $f$ gives rise to the multiplication of
the series $W$ by $f(q\exp(\h z)) f(q)$ which obviously preserves the
polynomiality property of the coefficients $w_{d,m}$. This proves (ii).

\medskip

{\bf Proposition 7.2.} {\em Let $f=\sum_{d\in \L'_0\cup \L_1} f_d q^d$ be
a series of weighted - homogeneous degree $1$ where $f_d$ 
are linear inhomogeneous functions of $\l,\l'$ with rational coefficients,
and $f_0=0$. Then $\exp (f/\h ) Z$ satisfies the conditions (i),(ii).}

{\em Proof.} We have 
\[ \exp (-nf/u_j(\a) - f/\h)= \exp (-(n\h+u_j(\a) f/u_j(\a) \h) =
1+ (n\h +u_j(\a)) g_{\a,j,n} ,\] 
where $g_{\a,j,n}$ is a $q$-series with coefficients
polynomial in $\h^{-1}$ and rational in $(\l,\l')$.  
This implies that $\exp (f/h ) Z_{\a}$ satisfy the recursion relation (i) with
the new initial condition 
\[ e^{f/\h} \sum_{d\in \L} R_{d,\a}q^d + \]
\[e^{f/\h} \sum_{j\notin \a}\sum_{n=1}^{\infty}
g_{\a, j, n} \ q^{n\ d_{\a,j}} \ C_{\a,j}(n) \ 
Z_{\b(\a,j)}(q,-u_j(\a)/n;\l,\l')  \]
whose coefficients at $q^d$ are still polynomial in $1/\h$. 

The transformation $Z\mapsto \exp (f/\h) Z$ multiplies the series $W$ by
$\exp g$, where $g=(f(qe^{\h z})-f(q))/\h$ has polynomial
coessicients when written as a power series in $(q,z)$ since $\exp (\h z d)-1$
is divisible by $\h$. $\square $

\medskip

{\bf Proposition 7.3.}{\em Let the linear combination $fp=f_1 p_1+...+f_k p_k$ 
of the equivariant cohomology classes $p_i\in H^*_T(Y,\QQ )$ be given by the
series $f_i=\sum_{d\in \L'_0} f_{i,d}q^d$  
with rational coefficients $f_{i,d}\in \QQ $ and $f_{i,0}=0$. Then
$ \exp (fp/\h) Z(q\exp f ,\h ;\l ,\l')$ satisfies the conditions (i),(ii).}

{\em Proof.} For $\b=\b(\a,j)$ we have $p(\a)-p(\b)=d_{\a,j} u_j(\a)$ and 
therefore
\[ \exp [-fp(\a)/\h - nfp(\b)/u_j(\a)] =\]
\[ e^{n\ d_{\a,j} f}\ \exp [-fp(\a) (n\h+u_j(\a))/(u_j(\a)\h )] = 
 e^{n\ d_{\a,j} f} + (n\h+u_j(\a)) g_{\a,j,n} \ , \]
where the series $g_{\a,j,n} (q,\h^{-1};\l,\l')$ has the same properties as
specified in the proof of Proposition $7.2$. This implies that the recursion
relation for $e^{fp/\h}Z(qe^f,\h)$ (with some initial conditions
$\sum R_{\a,d}q^d$) follows from the recursion relation
obtained by the change $q\mapsto qe^f$ from (i) satisfied by $Z$.
 
The operation $Z\mapsto e^{fp}Z(qe^f,\h)$ transforms $W(z,q)$ to
\[ W(z+\frac{f(qe^{\h z})-f(q)}{\h} , qe^f ) \]
and preserves the polynomiality property of the coefficients $w_{d,m}$
since coefficients of the series $f_i(qe^{\h z})-f_i(q)$ 
are divisible by $\h $. $\square $

\medskip

Using the notations from Proposition $6.4$ we perform the following 
operations with the series $\Psi $ which according to Propositions $7.1$,
$7.2$, $7.3$ preserve the properties (i),(ii).

(1) Divide $\Psi $ by $\Psi^{(0)}$; this transforms the 
asymptotical expansion $\Psi=\Psi^{(0)}+\Psi^{(1)}/\h + o (1/\h)$ to
\[ 1+ \h^{-1}\ (H(q) + \sum_a G_a(q) v_a -\sum_j F_j(q) u_j)
+ o (\h^{-1}), \]
where  
$H=(\sum_{d\in \L_1}\frac{L_1!...L_l!}{D_1!...D_N!}q^d)/\Psi^{(0)}$,
and $G_a=g_a/\Psi^{(0)}, \ F_j=f_j/\Psi^{(0)}$.
 
(2) Divide $\Psi /\Psi^{(0)}$ by 
$\exp \h^{-1} (H+\sum_j F_j\l_j-\sum_a G_a\l'_a ) $;
this transforms the asymptotical expansion to
\[ Z = 1+\h^{-1}\ \sum_i \phi_i p_i + o (\h^{-1}) \ ,\]
where $\phi_i=\sum _a l_{ia} G_a - \sum_j m_{ij} F_j$.

(3) Transform $\exp (-\sum \phi_ip_i/\h) Z (q,\h )$ to the new variables 
$Q_i=q_i e^{\phi_i(q)}$, $i = 1,...,k$; the resulting series $S(Q,\h;\l,\l')$
has the asymptotical expansion $1+o(1/h)$.

\medskip
  
{\bf Corollary 7.4.} {\em Suppose that $\calt_{Y}$ is non-negative. Then
$\cals (q,\h;\l,\l')=S(q,\h;\l,\l')$.}

{\em Proof.} The series $S$ satisfies the conditions (a),(b),(c) 
of Proposition $4.5$ which uniquely determine the series $\cals $.  


\enddocument